\input ./style/arxiv-general.cfg
\documentclass[aos,MSNbibl,nameyear,seceqn,dvips]{arximspdf}
\makeatletter
   \@ifpackageloaded{graphicx}{}{\usepackage{graphicx}}
\makeatother
\usepackage{dcolumn}
\usepackage{mathrsfs}
\doi{10.1214/15-AOS1372}% Updated by VTEXPTS2LaTeX.exe, 06.10.2015 12:51
\volume{44}
\issue{1}
\pubyear{2016}
\firstpage{425}
\lastpage{453}
\docsubty{FLA}

\makeatletter
\newcolumntype{d}[1]{D{.}{.}{#1}}
\newcommand{\rrVert}{\Vert}
\newcommand{\llVert}{\Vert}

\newcommand{\dddot}[1]{\stackrel{\cdot\cdot\cdot}{#1}}
\renewcommand{\citep}[1]{(\citeauthor{#1} \citeyear{#1})}
\newcommand{\eqref}[1]{(\ref{#1})}
\newproclaim{rem}{Remark}[section]
\newproclaim{defn}{Definition}[section]
\newtheorem{theo}{Theorem}[section]

\newtheorem{cor}{Corollary}[section]

\newtheorem{lemma}{Lemma}[section]
\newproclaim{ass}{Assumtpion}
\newproclaim{example}{Example}[section]
\makeatother

\begin{document}
\begin{frontmatter}

\title{Estimation for single-index and partially linear single-index
integrated models}
\runtitle{Partially linear single-index models}

\begin{aug}
% Corresponding author: Dag Tjostheim - Dag.Tjostheim@math.uib.no% Updated by VTEXPTS2LaTeX.exe, 07.10.2015 07:32
%Updated by VTEXPTS2LaTeX.exe, 06.10.2015 12:51
\author[A]{\fnms{Chaohua}~\snm{Dong}\thanksref{M1,M2,T1}\ead[label=e1]{dchaohua@hotmail.com}},
\author[B]{\fnms{Jiti}~\snm{Gao}\thanksref{M2,T1}\ead[label=e2]{jiti.gao@monash.edu}}
\and
\author[C]{\fnms{Dag}~\snm{Tj{\o}stheim}\corref{}\thanksref{M3}\ead[label=e3]{dag.tjostheim@math.uib.no}}
\runauthor{C. Dong, J. Gao and D. Tj{\o}stheim}
\affiliation{Southwestern University of Finance and
Economics,\thanksmark{M1} Monash University\thanksmark{M2}  and
University of Bergen\thanksmark{M3}}
\thankstext{T1}{Supported in part by two
Australian Research Council Discovery Grants: \#DP1096374 and \#DP130104229.}
%\runauthor{}
%\dedicated{}
\address[A]{C. Dong\\
School of Economics\\
Southwestern University\\
\quad of Finance and Economics\\
Chengdu\\
China\\
and\\
Department of Econometrics \\
\quad and Business Statistics\\
Monash University\\
Melbourne\\
Australia\\
\printead{e1}}

\address[B]{J. Gao\\
Department of Econometrics \\
\quad and Business Statistics\\
Monash University\\
Melbourne\\
Australia\\
\printead{e2}}

\address[C]{D. Tj{\o}stheim\\
Department of Mathematics\\
University of Bergen\\
5020 Bergen\\
Norway\\
\printead{e3}}
\end{aug}

% HISTORY:
%
\received{\smonth{1} \syear{2015}}% Updated by VTEXPTS2LaTeX.exe,
%06.10.2015 12:51
%
\revised{\smonth{8} \syear{2015}}% Updated by VTEXPTS2LaTeX.exe,
%06.10.2015 12:51

% ABSTRACT
%
\begin{abstract}
Estimation mainly for two classes of popular models, single-index and
partially linear single-index models, is studied in this paper. Such
models feature nonstationarity. Orthogonal series expansion is used to
approximate the unknown integrable link functions in the models and a
profile approach is used to derive the estimators. The findings include
the dual rate of convergence of the estimators for the single-index
models and a trio of convergence rates for the partially linear
single-index models. A new central limit theorem is established for a
plug-in estimator of the unknown link function. Meanwhile, a
considerable extension to a class of partially nonlinear single-index
models is discussed in Section~\ref{sec4}. Monte Carlo simulation verifies these
theoretical results. An empirical study furnishes an application of the
proposed estimation procedures in practice.
\end{abstract}

% KEYWORDS
%
%\begin{keyword}
%\kwd{}
%\end{keyword}
\end{frontmatter}

%s1 #&#
\section{Introduction}\label{sec1}
In the last decade or so, nonlinear (nonparametric or semiparametric)
and nonstationary time series models have been studied extensively and
improved dramatically as witnessed by the literature, such as those
based on the nonparametric kernel approach by \citet{karlsen2001,
karlsen2007}, \citeauthor{jiti2009a} (\citeyear{jiti2009b,jiti2009a}), \citet{phillips2009},
\citeauthor{qiyingwang2009a} (\citeyear{qiyingwang2009a,qiyingwang2009b,qiyingwang2012}), \citet{jiti2012a},
Gao and Phillips (\citeyear{gp13,jiti2013}) and \citet{jiti2013a}, among
others. The main development in the field is the establishment of new
estimation and specification testing procedures as well as the
resulting asymptotic properties. In recent years, the conventional
nonparametric kernel-based estimation and specification testing theory
has been extended to the nonparametric series based approach; see, for
example, \citeauthor{donggao2013} (\citeyear{donggao2013,donggao2014}).

We first consider a partially linear single-index model of the form
%
%e1.1 #&#
\begin{equation}
\label{a1} y_t=\beta_0^\top
x_t+g\bigl(\theta_0^\top x_t
\bigr)+e_t,\qquad t=1,\ldots, n,
\end{equation}
where $y_t$ is a scalar process, $g(\cdot)$, the so-called link
function, is an unknown nonlinear integrable function from $\mathbb
{R}$ to $\mathbb{R}$, $\beta_0$ and $\theta_0$ are the true but unknown
$d$-dimensional column vectors of parameters, the superscript $^\top$
signifies the transpose of a vector (or matrix, hereafter), $x_t$ is a
$d$-dimensional integrated process, $e_t$ is an error process and $n$
is sample size.

The motivations of this study are as follows. In a fully nonparametric
estimation context, researchers often suffer from the so-called ``curse
of dimensionality'', and hence dimensionality reduction is particularly
of importance in such a situation. One efficient way of doing so is to
use index models like model \eqref{a1}. Moreover, model \eqref{a1} is
also an extension of linear parametric models, since it would become a
linear model under the particular choice of the link function. Taking
these into account, models such as \eqref{a1} are often used as a
reasonable compromise between fully parametric and fully nonparametric
modelling. See, for example, \citet{fanjianqing1997}, \citet{xia1999,
xia2002}, \citet{yu2002}, \citet{zhu2006}, \citet{lianghua2010}, \citet
{wangjane2010} and \citet{ma2013}. Nevertheless, most researchers only
focus on the stationary covariate case so that their theoretical
results are not applicable for practitioners who use partially linear
single-index model to deal with nonstationary time series data. For
example, in macroeconomic context practitioners may be concerned with
inflation, unemployment rates and other economic indicators. These
variables exhibit nonstationary characteristics. Therefore, it is
desirable in such circumstances to develop estimation theory for the
partially linear single-index models.

Furthermore, recent studies by \citeauthor{gp13} (\citeyear{gp13,jiti2013}) have\break pointed out
that, for multivariate $I(1)$ processes, the conventional kernel
estimation method may not be workable because the limit theory may
break down. This gives rise to a challenge of seeking alternative
estimation methods.

When $\beta_0=0$, model \eqref{a1} becomes a single-index model
%
%e1.2 #&#
\begin{equation}
\label{a2} y_t=g\bigl(\theta_0^\top
x_t\bigr)+e_t,\qquad t=1,\ldots, n,
\end{equation}
which has been studied extensively for the case where $x_t$ is
stationary [see, e.g., \citet{ichimura1993}, \citet{xia1999b} and \citet
{yuyan2010}].

We shall first consider model \eqref{a2}, but this is mainly a
preliminary stage on our way to the general model \eqref{a1}. Standing
on its own, model \eqref{a2} has limited applicability since it is
integrable, and among other things, does not include the linear model.
Coupled with the assumption that $\{x_t\}$ is a unit root process, this
implies that only an order of $O(\sqrt{n})$ observations can be used in
the estimations of $\theta_0$ and $g$ in \eqref{a2}. The function $g$
is only capable of describing finite domain behaviour in $x_t$. As
$\theta_0^\top x_t$ increases, $g(\theta_0^\top x_t)$ goes to zero. All
of this will be made precise in the following.

When the $g$ function is added as a component in model \eqref{a1}, one
obtains a model whose behaviour is governed by the linear component
with $g$ superimposed, whereas as $x_t$ becomes large it reduces to the
linear component. The resulting model \eqref{a1} can be likened to a
smooth transition regression model [STR model, \citet{dag2010}], where
as $x_t$ increases the model changes smoothly to the linear model. Our
model in this sense extends the smooth transition model to a situation
where there is an index involved and with a nonstationary input
process. We believe that this is a situation which is of interest both
theoretically and practically, as witnessed for example by our
empirical study where quite different index behaviour is obtained in
the close domain as compared to the far out region.

It is indeed possible to generalise our model to include a nonlinear
behaviour also far out, which leads to the next stage of modelling. As
a linear function is a particular $H$-regular function while an
integrable function belongs to $I$-regular functions, studied by \citeauthor{phillips1999}
(\citeyear{phillips1999,phillips2001}), \citet{qiyingwang2009a}, in a third stage we
extend model \eqref{a1} using a known $H$-regular function to
substitute the linear function, in order to make the model more
flexible and applicable.

Following the existing identifiability condition, for example, \citet
{linwei2007}, we assume for models studied later that $\|\theta_0\|=1$
and the first nonzero component of $\theta_0$ is positive. Notice that
there is no extra condition needed for $\beta_0$ to make \eqref{a1}
identifiable, as discussed in Section~\ref{sec2.2} below. To facilitate the
theoretical development in the following sections, we assume that
$\theta_0$ is an interior point located within a compact and convex
parameter space $\Theta$, which is also a usual assumption in a
parameter estimation context. To focus on the unit root case, we also
assume throughout that cointegration will not happen for $\theta$
around $\theta_0$. In other words, there exists a neighbourhood of
$\theta_0$, $\mathcal{N}(\theta_0, \delta)\subset\Theta$, such that
for any $\theta\in\mathcal{N}(\theta_0, \delta)$, $\theta^\top x_t$ is
always an $I(1)$ process.

The findings of this paper are summarised as follows. The rate of
convergence of the estimators of $\theta_0$ in both models \eqref{a1}
and \eqref{a2} is a composite of two different rates in a new
coordinate system where $\theta_0$ is on one axis. $\widehat{\theta}_n$
has a rate of $n^{-1/4}$ on the $\theta_0$-axis, and another rate as
fast as $n^{-3/4}$ on all axes orthogonal to $\theta_0$. Overall,
$\widehat{\theta}_n$ possesses convergence rate $n^{-1/4}$. This is
expected and comes from the integrability of $g(\cdot)$, which in turn
reduces the number of effective observations to $\sqrt{n}$. Moreover,
the rate of convergence of $\widehat{\beta}_n$ to $\beta_0$ is
$n^{-1}$, consistent with that of a linear model with a unit root
input. The normalisation of $\widehat{\theta}_n$, $\|\widehat{\theta
}_n\|^{-1}\widehat{\theta}_n$, converges to $\theta_0$ with a rate
faster than $\widehat{\theta}_n$ in both models. A new central limit
theorem for a plug-in estimator of the form $\widehat{g}_{n} (u)$
converging to $g(u)$, where $u\in\mathbb{R}$, is comparable with the
conventional kernel estimator in the literature. These phenomena are
verified with finite sample experiments below.

Theoretical results heavily depend on the level of nonstationarity of
the integrated time series and the integrability of the link functions.
These properties result in a slow rate of convergence for the link
function involving an $I(1)$ process and fast rate of convergence for a
linear model with an $I(1)$ process. These are very different from the
literature where the regressors are stationary. In addition, Monte
Carlo simulations generally need relative larger sample sizes than
those for the cases where the regressors are stationary if the
regression function is integrable, since random walk on one hand
diverges at rate $\sqrt{n,}$ and on the other hand it possesses
recurrent property making it possible to return to the effective domain
of the integrable function $g$.

Two papers related to this study are \citet{chang2003} and \citet
{moon2006} in terms of regressor. However, \citet{ chang2003} stipulate
that their link function is a smooth distribution function-like
transformation and they are not interested in the estimation of the
unknown link function. \citet{moon2006} point out in the discussion
section that their method developed for binary choice models may be
applicable for the estimation of single-index models where the link
function $g(x)\to\infty$ as $|x|\to\infty$. Clearly, they are quite
different from the setting of this study.

The organisation of the rest of the paper is as follows. Section~\ref{sec2}
gives estimation procedures and assumptions for models \eqref{a1} and
\eqref{a2}. Asymptotic theory is established in Section~\ref{sec3} for the
estimator $\widehat{\theta}_n$ in model \eqref{a2} and the estimator
($\widehat{\beta}_n$, $\widehat{\theta}_n$) in model \eqref{a1}. A
central limit theorem for a plug-in estimator of the form $\widehat
{g}_{n} (u)$ is given in Section~\ref{sec3}. An extension of model \eqref{a1} is
discussed in Section~\ref{sec4} and Monte Carlo simulation experiments are
conducted in Section~\ref{sec5}. Section~\ref{sec6} shows the implementation of the
proposed estimation schedules with an empirical dataset. Appendix \ref{appa}
presents some technical lemmas. The proof of the main results in
Section~\ref{sec3} is given in Appendix \ref{appb}. A supplemental document [\citet
{DGT2015}] contains Appendices C, D and E where all the proofs of the
key lemmas listed in Appendix \ref{appa} as well as some other lemmas are shown
in Appendix C, the complete proof of the results in Section~\ref{sec3} is placed
in Appendix~D and the results in Section~\ref{sec4} are proven in Appendix E.

Throughout the paper, the following notation is used. $\|\cdot\|$ is
Euclidean norm for vectors and element-wise norm for matrices, that is,
if $A= (a_{ij})_{nm}$, $\|A\|=(\sum_{i=1}^n \sum_{j=1}^m a_{ij}^2
)^{1/2}$; $I_d$ is the $d$-dimensional identity matrix; $[a]$ is the
maximum integer not exceeding $a$; $\mathbb{R}$ is the real line; for
any function $f(\cdot)$, $\dot{f}(x)$, $\ddot{f}(x)$ and $\dddot{f}(x)$
are the derivatives of the first, second and third order of $f(\cdot)$
at $x$. Here, when $f(x)$ is a vector-valued function its derivatives
should be understood as element-wise. Furthermore, $\phi(\cdot)$ stands
for the density function of a multivariate standard normal variable;
$\int f(w)\,dw$ means a multiple integral when $w$ is a vector.
Convergence in probability and convergence in distribution are
signified as $\to_P$ and $\to_D$, respectively.

%s2 #&#
\section{Estimation procedure and assumptions}\label{sec2}

Suppose that the link function $g(\cdot)$ belongs to $L^2(\mathbb{R}) =
\{f(x): \int f^2(x)\,dx <\infty\}$. It is known that the Hermite
function sequence $\{\mathscr{H}_i(x)\}$ is an orthonormal basis in
$L^2(\mathbb{R})$ where by definition
%
%e2.1 #&#
\begin{equation}
\label{b1} \mathscr{H}_i(x)=\bigl(\sqrt{\pi} 2^ii!
\bigr)^{-1/2} H_i(x)\exp \biggl(-\frac
{x^2}{2} \biggr), \qquad i
\ge0,
\end{equation}
and $H_i(x)$ are Hermite polynomials orthogonal with density $\exp
(-x^2)$. The orthogonality reads $\int\mathscr{H}_i(x)\mathscr
{H}_j(x)\,dx=\delta_{ij}$, the Kronecker delta.

Thus, a continuous function $g(\cdot)\in L^2(\mathbb{R})$ may be
expanded into an orthogonal series
%
%e2.2 #&#
\begin{equation}
\label{b3} g(x)=\sum_{i=0}^\infty
c_i\mathscr{H}_i(x) \quad\mbox{and}\quad c_i=\int g(x)
\mathscr{H}_i(x)\,dx.
\end{equation}
Throughout, let $k$ be a positive integer and define $g_k(x)=\sum_{i=0}^{k-1}c_{i}\mathscr{H}_i(x)$ as the truncation series and $\gamma
_k(x)=g(x)-g_k(x)=\sum_{i=k}^{\infty} c_i\mathscr{H}_i (x)$ as the
residue after truncation.

%s2.1 #&#
\subsection{Estimation procedure for single-index models}\label{sec2.1}

By virtue of \eqref{b3}, we write model \eqref{a2} for $t=1, \ldots, n$ as
\[
y_t=Z_k^\top\bigl(\theta^\top_0
x_t\bigr) c+\gamma_k\bigl(\theta^\top_0
x_t\bigr)+e_t,
\]
where $Z_k^\top(\cdot)= (\mathscr{H}_0(\cdot), \ldots, \mathscr
{H}_{k-1} (\cdot))$, $c^\top=(c_0, \ldots, c_{k-1})$ and $k$ is the
truncation parameter determined later.

Let $Y=(y_1,\ldots, y_n)^\top$, $Z=(Z_k(\theta^\top_0 x_1), \ldots,
Z_k(\theta^\top_0 x_n) )^\top$ an $n\times k$ matrix, $\gamma=(\gamma
_k(\theta^\top_0 x_1), \ldots, \gamma_k (\theta^\top_0 x_n))^\top$ and
$e=(e_1, \ldots, e_n)^\top$. We have a matrix form equation $Y=Zc+\gamma
+e$, and hence by the Ordinary Least Squares (OLS) method, $\widetilde
{c}= \widetilde{c} (\theta_0)=(Z^\top Z)^{-1}Z^\top Y$ is an estimate
for $c$ in terms of $\theta_0$. Nonetheless, since $\theta_0$ is
unknown, we only have a form of $\widetilde{c}$. To estimate $\theta
_0$, define for $\theta\in\Theta$, $L_n(\theta)= \frac{1}{2}\sum_{t=1}^n [y_t-Z_k^\top(\theta^\top x_t) \widetilde{c}(\theta)]^2$.
Then we choose an optimum $\widehat{\theta}_n$ such that
%
%e2.3 #&#
\begin{equation}
\label{b5} \widehat{\theta}_n= \mathop{\operatorname{argmin}}_{ \theta\in\Theta}
L_n(\theta),
\end{equation}
as an estimator for $\theta_0$. Once $\widehat{\theta}_n$ is available,
we have a plug-in estimator $\widehat{g}_n(u) \equiv\widehat{g}_n(u;
\widehat{\theta}_n)=Z_k(u)^\top\widehat{c}$ for any $u\in\mathbb{R}$
where $\widehat{c}= \widetilde{c}(\widehat{\theta}_n)$, which is purely
based on the sample, and hence applicable. The estimation procedure
proposed here is the profile method [see, \citet{severini1992, lianghua2010}].

Additionally, to be in concert with the identification condition $\|
\theta_0\|=1$, we define the normalisation of $\widehat{\theta}_n$,
$\widehat{\theta}_{n,\mathrm{emp}}=\|\widehat{\theta}_n\|^{-1} \widehat
{\theta}_n$. An asymptotic theory for both $\widehat{\theta}_{n,\mathrm
{emp}}$ and $\widehat{\theta}_n$ will be studied in Section~\ref{sec3} below.

%s2.2 #&#
\subsection{Estimation procedure for partially linear single-index models}\label{sec2.2}

Usually, researchers, such as \citet{xia1999}, impose an identification
condition that $\beta_0$ is perpendicular to $\theta_0$ on the
partially linear single-index models. This is because when $\beta_0$ is
not perpendicular to $\theta_0$, a new vector $\beta_0- (\beta_0^\top
\theta_0)\theta_0$ can be used in the place of $\beta_0$ and the $g$
function will be replaced by $g(u)+ (\beta_0^\top\theta_0) u$.
However, in model \eqref{a1} the lack of orthogonality between $\beta
_0$ and $\theta_0$ does not affect the identifiability of the model at
all. See the verification at the end of Appendix C in the supplementary
material [\citet{DGT2015}].

Our estimation procedure in partially linear single-index models is
proposed as follows. By virtue of \eqref{b3} again, for each $t$
rewrite \eqref{a1} as
\[
y_t-\beta^\top_0 x_t=Z_k
\bigl(\theta^\top_0 x_t\bigr)^\top c+
\gamma_k\bigl(\theta^\top _0 x_t
\bigr)+e_t,
\]
where $Z_k(\cdot)$, $c$ and $\gamma_k(\cdot)$ are defined as before.

Denote $X=(x_1, x_2, \ldots, x_n)^\top$ an $n\times d$ matrix, and $Y,
Z, \gamma, e$ remain the same as in the last subsection. We have matrix
form equation: $Y-X\beta=Zc+\gamma+ e$. Then the OLS gives that
$\widetilde{c}= \widetilde{c}(\beta_0,\theta_0)=(Z^\top Z)^{-1}Z^\top
(Y-X\beta_0)$. Due to the same reason as before, define for generic
$(\beta, \theta)$, $L_n(\beta, \theta)=\break\frac{1}{2} \sum_{t=1}^n
[y_t-\beta^\top x_t-Z_k^\top(\theta^\top x_t) \widetilde{c}(\beta
,\theta)]^2$. The estimator of $(\beta_0, \theta_0)$ is given by
%
%e2.4 #&#
\begin{equation}
\label{b9} %
\pmatrix{\widehat{\beta}_n \vspace*{2pt}
\cr
\widehat{\theta}_n } %
=\mathop{\operatorname{argmin}}_{\beta\in\mathbb{R}^d, \theta\in\Theta}
L_n(\beta, \theta).
\end{equation}
Similarly, a plug-in estimator is obtained, $\widehat{g}_n(u)\equiv
\widehat{g}_n(u; \widehat{\beta}_n, \widehat{\theta}_n) =Z_k^\top
(u)\widehat{c}$ where $\widehat{c}= \widetilde{c}(\widehat{\beta}_n,
\widehat{\theta}_n)$. Once the estimators of the parameters are
available, and the \mbox{normalisation} $\widehat{\theta}_{n,\mathrm{emp}}=\|
\widehat{\theta}_n\|^{-1} \widehat{\theta}_n$ is defined to satisfy the
identification condition.

%s2.3 #&#
\subsection{Assumptions}\label{sec2.3}
Before we establish our main theory in Section~\ref{sec3} below, we introduce
some necessary conditions.

\renewcommand{\theass}{\Alph{ass}}
\setcounter{ass}{0}
%as1 #&#
\begin{ass}\label{assa}
%\textbf{Assumption A}
%%
(a) Let $\{\varepsilon_j, -\infty<j< \infty\}$ be a sequence of
$d$-dimen\-sional independent and identically distributed (i.i.d.)
continuous random variables with $E\varepsilon_1=0$, $E[\varepsilon_1
\varepsilon_1^\top] =\Omega>0$ and $E\|\varepsilon_1 \|^p <\infty$ for
some $p>2$. The characteristic function of $\varepsilon_1$ is
integrable, that is, $\int|E\exp(iu \varepsilon_1)| \,du <\infty$.\vspace*{-6pt}
\begin{longlist}[(b)]
\item[(b)] Let $x_t=x_{t-1}+v_t$ for $t\geq1$ and $x_0=O_P(1)$, where $\{
v_t\}$ is a linear process defined by $v_t=\sum_{j=0}^\infty\rho_j
\varepsilon_{t-j}$, in which $\{\rho_j\}$ is a square matrix such that
$\rho_0=I_d$, $\sum_{j=0}^\infty\|\rho_j\|<\infty$ and $\rho= \sum_{j=0}^\infty\rho_j$ is of full rank.
\item[(c)] There is a $\sigma$-field $\mathcal{F}_{t}$ such that $(e_t,
\mathcal{F}_{t})$ is a martingale difference sequence, that is, for all
$t$, $E(e_t|\mathcal{F}_{t-1}) =0$ almost surely (a.s.). Also,\break
$E(e^2_t| \mathcal{F}_{t-1}) =\sigma^2_e$ a.s. and $\mu_4:= \sup_{1\le
t\le n} E(e_t^4| \mathcal{F}_{t-1})<\infty$ a.s.
\item[(d)]$x_{t}$ is adapted with $\mathcal{F}_{t-1}$.
\item[(e)] Let $V_n(r)=\frac{1}{\sqrt{n}}\sum_{i=1}^{[nr]}v_i$ and
$U_n(r)=\frac{1}{\sqrt{n}} \sum_{i=1}^{[nr]} e_i$. Suppose that
$(U_n(r),\break  V_n(r))\rightarrow_D (U(r),V(r))$ as $n\to\infty$. Here,
$(U(r),V(r))$ is a $(d+1)$-vector of Brownian motions.
\end{longlist}
\end{ass}

%re2.1 #&#
\begin{rem}
All conditions in Assumption \ref{assa} are routine requirements in the
nonstationary model estimation context. Conditions (a) and (b)
stipulate that the regressor $x_t$ is an integrated process generated
by a linear process $v_t$ which has the i.i.d. sequence $\{\varepsilon
_j, -\infty<j< \infty\}$ as building blocks. Meanwhile, (c), (d) and
(e) are extensively used in related papers such as \citet
{phillips2000}, \citeauthor{qiyingwang2009a} (\citeyear{qiyingwang2009a,qiyingwang2009b,qiyingwang2012}),
\citeauthor{jiti2009b} (\citeyear{jiti2009b,jiti2009a}), \citet{jiti2012}, among others.
The $\sigma$-field $\mathcal{F}_{t}$ may be taken as $\mathcal
{F}_{t}=\sigma( \ldots, \varepsilon_{t}, \varepsilon_{t+1}; e_1,\ldots, e_t)$.
\end{rem}

By Skorohod representation theorem [\citet{pollard1984}, page 71] there
exists $(U^0_n(r), V_n^0(r))$ in a richer probability space such that
$(U_n(r), V_n(r)) =_D (U^0_n(r), V_n^0(r))$ for which $(U^0_n(r),
V_n^0(r))\to_{a.s.} (U(r), V(r))$ uniformly on\break $[0,1]^{d+1}$. To avoid
the repetitious embedding procedure of $(U_n(r), V_n(r))$ to the richer
probability space where $(U^0_n(r), V_n^0(r))$ is defined, we simply
write $(U_n(r),\break V_n(r))= (U^0_n(r), V_n^0(r))$ instead of $(U_n(r),
V_n(r))=_D (U^0_n(r), V_n^0(r))$. Since Lemma~\ref{lemma2} below is
derived in this richer probability space, all proofs in the paper
should be understood in the richer space as well. We will not repeat
this again.

%as2 #&#
\begin{ass}\label{assb}
%\textbf{Assumption B}
(a) $g(x)$ is differentiable on $\mathbb{R}$ and $g^{(m-\ell
)}(x)x^\ell\in L^2(\mathbb{R})$ for $\ell=0,1,\ldots, m$ with some
given integer $m$.\vspace*{-6pt}
\begin{longlist}[(b)]
\item[(b)]$k=[a\cdot n^\kappa]$ with some constant $a>0$, $\kappa\in
(0,1/8)$ and $\kappa(m-3)\ge\frac{1}{2}$ with $m$ as in (a) above.
\end{longlist}
\end{ass}

%re2.2 #&#
\begin{rem}
Condition (a) ensures the negligibility of the truncation residuals
(see the derivation at the beginning of Lemma C.1 of Appendix C of the
supplementary document). Regarding condition (b), although it is
stringent for $\kappa$, we may choose, for example, $\kappa\in[\frac
{5}{44}, \frac{5}{41}]$ and $m=8$ in practice. Large $m$ and small
$\kappa$ are chosen such that the orthogonal series expansion for the
link function converges so fast that all residues after truncation do
not affect the limit theory, as can be seen in the proof of Theorem~\ref{th31}
in the supplementary
material \citet{DGT2015}.
\end{rem}

%s3 #&#
\section{Asymptotic theory}\label{sec3}

%s3.1 #&#
\subsection{Asymptotic theory for single-index models}\label{sec3.1}

To derive an asymptotic theory for $\widehat{\theta}_n$ given by \eqref
{b5}, we shall use basic ideas from \citet{wooldridge1994}. Let
$S_n(\theta)=\frac{\partial}{\partial\theta}L_n(\theta)$ and
$J_n(\theta)= \frac{\partial^2 }{\partial\theta\,\partial\theta^\top}
L_n(\theta)$ be the score and Hessian, respectively. As usual, we have
the expansion
%
%e3.1 #&#
\begin{equation}
\label{3a} 0=S_n(\widehat{\theta}_n)=S_n(
\theta_0)+J_n(\theta_n) (\widehat{\theta
}_n -\theta_0),
\end{equation}
where $J_n(\theta_n)$ is the Hessian matrix with the rows evaluated at
a point $\theta_n$ between $\widehat{\theta}_n$ and $\theta_0$.

To facilitate the development of the asymptotic theory, we consider
coordinate rotation in $\mathbb{R}^d$. Let $Q=(\theta_{0}, Q_2)$ be an
orthogonal matrix. Note that such $Q$ does exist since $\theta_0\ne0$.
We shall use the orthogonal matrix $Q$ to rotate all vectors in~$\mathbb
{R}^d$. In particular,
%
%e3.2 #&#
\begin{eqnarray}
\label{3b} %
\alpha_0&:=& Q^\top
\theta_0=\bigl(\alpha_{10}, \alpha_{20}^\top
\bigr)^\top \qquad\mbox{where } \alpha_{10}=\|\theta_{0}
\|^2=1, \alpha_{20}=Q_2^\top
\theta_{0}=0,\nonumber
\\
\qquad z_t&:=&Q^\top x_t=\bigl(x_{1t},
x_{2t}^\top\bigr)^\top\qquad \mbox{where }
x_{1t}:=\theta_0^\top x_t,
x_{2t}:=Q_2^\top x_t,
\\
\alpha&:=&Q^\top\theta \qquad\mbox{for any generic } \theta.\nonumber %
\end{eqnarray}

Accordingly, we can rewrite the single-index model as $y_t=g(\theta
_0^\top QQ^\top x_t)+e_t= g(\alpha_{0}^\top z_{t}) +e_t$. In addition,
by Assumption \ref{assa} and the continuous mapping theorem, we have for $r\in[0,1]$,
%
%e3.3 #&#
\begin{equation}\qquad
\label{3c} \frac{1}{\sqrt{n}}x_{1[nr]}\to_D
V_1(r)=\theta_0^\top V(r)\quad \mbox {and}\quad
\frac{1}{\sqrt{n}}x_{2[nr]} \to_DV_2(r)=Q_2^\top
V(r).
\end{equation}

It is noteworthy that the rotation is not necessary in practice, as
shown in the simulation section, and it is also logically impossible
since $\theta_0$ is unknown. The rotation is only used as a tool to
derive an asymptotic theory for the proposed estimator.

If $\widehat{\alpha}_n$ is the nonlinear least squares estimator of
$\alpha_0$, then $\widehat{\alpha}_n =Q^\top\widehat{\theta}_n$.
Moreover, the score function $S_n(\alpha)$ and the Hessian $J_n(\alpha
)$ for the parameter $\alpha$ can be obtained from those for $\theta$.
More precisely, $S_n(\alpha)=Q^\top S_n(\theta)$ and $J_n(\alpha
)=Q^\top J_n( \theta)Q$. Premultiplying equation \eqref{3a} by $Q^\top
$, we have
%
%e3.4 #&#
\begin{equation}
\label{3d} 0=S_n(\widehat{\alpha}_n)=S_n(
\alpha_0)+J_n(\alpha_n) (\widehat{\alpha
}_n -\alpha_0).
\end{equation}

The following theorem gives asymptotic distributions for the score
$S_n(\alpha_0)$ and the Hessian $J_n(\alpha_0)$ as well as $\widehat
{\alpha}_n -\alpha_0$.

%th3.1 #&#
\begin{theo}\label{th31}
Denote $D_n= \operatorname{diag} (n^{1/4}, n^{3/4}I_{d-1})$. Under Assumptions
\ref{assa} and \ref{assb}, as $n\rightarrow\infty$
%
%e3.5 #&#
\begin{equation}
\label{th31a} D_n^{-1}S_n(
\alpha_0)\to_D R^{1/2}W(1) \quad\mbox{and}\quad
D_n^{-1}J_n(\alpha_0)D_n^{-1}
\to_P R,
\end{equation}
where $W(1)$ is a $d$-dimensional vector of standard normal random
variables independent of $V(r)$, and the symmetric block matrix $R= \bigl(
{r_{11}\atop r_{21}} \enskip {r_{12}\atop r_{22}}
 \bigr)$ is given by
\begin{eqnarray*}
r_{11}&=& L_1(1,0)\int s^2\dot{g}^2(s)\,ds,\qquad
r_{12}=\int_0^1V_2^\top
(r)\,dL_1(r,0)\int s\dot{g}^2(s)\,ds,
\\
r_{21}&=& r_{12}^\top,\qquad  r_{22}=\int
_0^1V_2(r)V_2^\top(r)\,dL_1(r,0)
\int \dot{g}^2(s)\,ds,
\end{eqnarray*}
in which $V_1$ and $V_2$ given by \eqref{3c} are Brownian motions of
dimension $1$ and $d-1$, respectively, $L_1(r,0)$ denotes the local
time process of Brownian motion $V_1(\cdot)$, standing for the
sojourning time of $V_1$ at zero over $[0,r]$.

As a result, under the same conditions, $\widehat{\alpha}_n$ is
consistent and as $n\to\infty$
%
%e3.6 #&#
\begin{equation}
\label{th31b} D_n(\widehat{\alpha}_n-
\alpha_0)\to_D R^{-1/2}W(1).
\end{equation}
\end{theo}

A standard book introducing the local time process of Brownian motion
is \citet{yor2005}. In view of the structure of $D_n$, we have two
limits from~\eqref{th31b},
%
%e3.7 #&#
\begin{equation}
\label{3e} n^{1/4}(\widehat{\alpha}_{1n}-1)
\to_D \mathbf{MN}(0,\rho_{11})\quad \mbox{and}\quad n^{3/4}
\widehat{\alpha}_{2n} \to_D \mathbf{MN}(0,
\rho_{22}),
\end{equation}
where $\widehat{\alpha}_{n}=(\widehat{\alpha}_{1n}, \widehat{\alpha
}_{2n}^\top)^\top$, ${\mathbf{MN}}(0, \Xi)$ stands for mixture normal
distribution for the case where the covariance matrix $\Xi$ is
stochastic, $\rho_{11}$ and $\rho_{22}$ are diagonal blocks on the
matrix $R^{-1}= \bigl(
{\rho_{11}\atop \rho_{21} }\enskip {\rho_{12}\atop \rho_{22}}
 \bigr)$,
%
%e3.8 #&#
\begin{equation}
\label{3f} \rho_{11}=\bigl(r_{11}-r_{12}r^{-1}_{22}r_{21}
\bigr)^{-1}\quad \mbox{and}\quad \rho_{22}=\bigl(r_{22}-r_{21}r_{11}^{-1}r_{12}
\bigr)^{-1}.
\end{equation}
Hence, $\widehat{\alpha}_{n}$ has two different convergence rates for
its components.

Note by \eqref{3e} that in the coordinate system $Q$ where $\theta_0$
is an axis, the estimator $\widehat{\theta}_n$ has dual convergence
rates: the rate of convergence for the coordinate on $\theta_0$ (i.e.,
$\widehat{\alpha}_{1n}$) is $n^{-1/4}$, while on all directions
orthogonal to $\theta_0$ the rate of convergence for the coordinates
(i.e., $\widehat{\alpha}_{2n}$) is as fast as $n^{-3/4}$. This
difference in convergence rate can be explained in the following way.
Due to the unit root behaviour of $\{x_t\}$ its probability mass is
spreading out in a Lebesgue type fashion. Since $g$ is integrable,
$g(\theta_0^\top x)\approx0$ outside the effective range of $g$. This
means that only moderate values of $\{x_t\}$ can contribute to $g$
along $\theta_0$, but in such directions that are orthogonal to $\theta
_0$, there is no restriction on $\{x_t\}$, so that even for far out
values of $\{x_t\}$, they can contribute, and hence increase the
effective sample size. Certainly, no such effect can take place in
univariate models.

As defined in Section~\ref{sec2}, $\widehat{\theta}_{n, \mathrm{emp}}= \|\widehat
{\theta}_n\|^{-1} \widehat{\theta}_n$. Intuitively, $\widehat{\theta
}_{n, \mathrm{emp}}$ might have a faster rate of convergence than that of
$\widehat{\theta}_n$. This can be seen using the $\alpha
$-representation of the rotated system. Because of $\widehat{\theta
}_n=Q \widehat{\alpha}_{n}$ and hence $\|\widehat{\theta}_n\|=\|
\widehat{\alpha}_{n}\|$, $\widehat{\theta}_{n, \mathrm{emp}}= Q\widehat
{\alpha}_{n, \mathrm{unit}}$, where $\widehat{\alpha}_{n, \mathrm{unit}}=(
\widehat{\alpha}_{n, \mathrm{unit}}^1, (\widehat{\alpha}_{n, \mathrm
{unit}}^2 )^\top)^\top=\|\widehat{\alpha}_{n}\|^{-1} \widehat{\alpha
}_{n}$. The following results give the rates of convergence for
$\widehat{\alpha}_{n, \mathrm{unit}}$ and then for $\widehat{\theta}_{n}$
and $\widehat{\theta}_{n, \mathrm{emp}}$, respectively.

%co3.1 #&#
\begin{cor}\label{cor1}
Under Assumptions \ref{assa} and \ref{assb}, we have as $n\to\infty$,
\[
n^{3/2} \bigl(\widehat{\alpha}_{n, \mathrm{unit}}^1-1\bigr)
\to_D -\tfrac{1}{2}\|\xi\| ^2 \quad\mbox{and}\quad
n^{3/4} \widehat{\alpha}_{n, \mathrm{unit}}^2\to_D
\xi,
\]
where $\xi\sim\mathbf{MN}(0, \rho_{22})$ is the limit given by
\eqref{3e}.
\end{cor}

Note that, after the normalisation, the slow rate becomes as fast as
$n^{-3/2}$ whereas the fast rate remains the same. Note also that
$\widehat{\alpha}_{n, \mathrm{unit}}^1\to_P 1$ but $\widehat{\alpha}_{n,
\mathrm{unit}}^1 =\|\widehat{\alpha}_{n} \|^{-1}\widehat{\alpha}_{1n}\le
1$. The intuitive reason for the fast rate of $\widehat{\alpha}_{n,
\mathrm{unit}}^1$ is that it takes advantage of the direction orthogonal
to $\theta_0$, where there is larger supply of information from far out
$x_t$'s as explained above. The rates are also verified with Monte
Carlo simulation. The resulting rates for $\widehat{\theta}_{n}$ are
given in Theorem~\ref{th32} below.

%th3.2 #&#
\begin{theo}\label{th32}
Under Assumptions \ref{assa} and \ref{assb}, we have as $n\to\infty$,
%
%e3.9 #&#
%e3.10 #&#
\begin{eqnarray}
n^{1/4} (\widehat{\theta}_n-\theta_0)
&\to_D&\mathbf{MN}\bigl(0,\rho_{11} \theta_0
\theta_0^\top\bigr), \label{th32a}
\\
n^{3/4} (\widehat{\theta}_{n,\mathrm{emp}}-\theta_0)
&\to_D& \mathbf {MN}\bigl(0, Q_2 \rho_{22}
Q_2^\top\bigr). \label{th32b}
\end{eqnarray}
\end{theo}

%re3.1 #&#
\begin{rem}
As can be seen, $\widehat{\theta}_n\to_P\theta_0$ at rate of $n^{-1/4}$
and $\widehat{\theta}_{n, \mathrm{emp}}\to_P \theta_0$ at rate of
$n^{-3/4}$. Again roughly speaking, the normalisation scales $\widehat
{\theta}_{n}$ to the unit ball, and hence accelerates the slow rate of
$\widehat{\theta}_{n}$. Due to the fast convergence of $\widehat{\theta
}_{n, \mathrm{emp}}$, all the following assertions regarding $\widehat
{\theta}_n$ remain true if $\widehat{\theta}_n$ is replaced by $\widehat
{\theta}_{n, \mathrm{emp}}$. A geometric illustration is given in
Appendix C of the supplementary material [\citet{DGT2015}] to explain
the slow and fast rates. We do not wish to repeat this again.
\end{rem}

Furthermore, by Theorem~\ref{th32} we have $\widehat{\theta}_n\sim\mathbf{MN}
(\theta_0, n^{-1/4} \rho_{11} \theta_0 \theta_0^\top)$. We
next show that the estimator of the covariance matrix of $\widehat
{\theta}_n$ is the inverse of the Hessian matrix of the form
$[J_n(\widehat{\theta}_n) ]^{-1}$ or even $[\widetilde{J}_n (\widehat
{\theta}_n) ]^{-1}$, where $\widetilde{J}_n(\theta)=\sum_{t=1}^n \dot
{\widehat{g}}{}^2_n (\theta^\top x_t)x_tx_t^\top$ is the leading term of
$J_n(\theta)$. Meanwhile, define the estimators for $\sigma_e$ and
$L_1(1,0)$ by
%
%e3.11 #&#
\begin{equation}\qquad
\label{3g} \widehat{\sigma}_e^2=\frac{1}{n}
\sum_{t=1}^n\bigl[y_t-
\widehat{g}_n\bigl(\widehat {\theta}_n^\top
x_t\bigr)\bigr]^2 \quad\mbox{and} \quad \widehat{L}_{n1}(1,0)=
\frac{1}{\sqrt{n}}\sum_{t=1}^n
\mathscr{H}_0^2 \bigl(\widehat{\theta}_n^\top
x_t\bigr),
\end{equation}
respectively, where $\mathscr{H}_0(\cdot)$ is the first function in the
Hermite sequence.

%co3.2 #&#
\begin{cor}\label{cor2}
Under Assumptions \ref{assa} and \ref{assb}, we have as $n\to\infty$,
%
%e3.12 #&#
\begin{equation}
\label{ehat} \widehat{\sigma}{}^2_e\to_P
\sigma_e^2\quad\mbox{and}\quad \widehat {L}_{n1}(1,0)-L_1(1,0)
\to_P 0
\end{equation}
and
%
%e3.13 #&#
\begin{equation}
\label{cor2a} \sqrt{n}\bigl[J_n(\widehat{\theta}_n)
\bigr]^{-1} \to_P \rho_{11} \theta_0
\theta _0^\top \quad\mbox{and}\quad \sqrt{n}\bigl[
\widetilde{J}_n (\widehat{\theta }_n)
\bigr]^{-1} \to_P \rho_{11} \theta_0
\theta_0^\top.
\end{equation}
\end{cor}

We then establish the following central limit theory for the plug-in
estimator $\widehat{g}_n(u)=Z_k^\top(u)\widehat{c}$ defined in Section~\ref{sec2.1}, where $u\in\mathbb{R}$.

%th3.3 #&#
\begin{theo}\label{th33}
Under Assumptions \ref{assa} and \ref{assb}, as $n\to\infty$, $\sup_{u\in\mathbb
{R}}|\widehat{g}_n(u)-g(u)|\to_P0$, and
%
%e3.14 #&#
\begin{equation}
\label{th33a} \widehat{\sigma}_e^{-1}
\widehat{L}_{n1}^{1/2}(1,0)n^{1/4}
\bigl\|Z_k(u)\bigr\| ^{-1}\bigl(\widehat{g}_n(u)-g(u)
\bigr)\to_DN(0, 1).
\end{equation}
\end{theo}

%re3.2 #&#
\begin{rem}
The order involved in the normality is $O_P(1)n^{1/4}k^{-1/2}$ in view
of $\|Z_k(u)\|^2= O(1)k$. This is comparable with the kernel estimate
in the literature. Theorem~3.1 of \citeauthor{qiyingwang2009a} [(\citeyear{qiyingwang2009a}), page~721]
shows that, for univariate regression $y_t=f(x_t)+u_t$, the normaliser
of $\hat{f}(x)-f(x)$ is $ (h\sum_{t=1}^n K_h(x_t-x)  )^{1/2}$
where $h$ is a bandwidth, $K_h(\cdot)=K(\cdot/h)/h$ is a kernel
function and $\hat{f}(x)$ is the kernel estimate of $f(x)$. Note that
$ (h\sum_{t=1}^n K_h(x_t-x) )^{1/2}=O_P(1)n^{1/4}h^{1/2}$.
Thinking of $k^{-1}$ as equivalent to the bandwidth $h$, the
normalisers in the two situations are quite comparable.
\end{rem}

%re3.3 #&#
\begin{rem}
Noting that $\widehat{g}_n(u)-g(u)=Z_k^\top(u)(\widehat{c}-c)-\gamma
_k(u)$ and by the orthogonality of the basis functions, $\int(\widehat
{g}_n(u)-g(u))^2\,dx= \|\widehat{c}-c\|^2+\|\gamma_k(u)\|^2$ where $\|
\gamma_k(u)\|$ is the norm of $\gamma(u)$ in the function space. Using
Lemma~\ref{cghat} below, $\|\widehat{c}-c\|^2=O_P(kn^{-1/2})$ and Lemma
C.1 (in the supplementary material), $\|\gamma_k(u)\|^2=O(k^{-m})$, an
optimal truncation parameter $k^*$ may be found to be proportional to
$k^*= [n^{1/2(m+1)} ]$ when $\|\widehat{c}-c\|^2$ and $\|\gamma
_k(u) \|^2$ have the same order going to zero. Here, $m$ is the
smoothness order of $g(u)$.
\end{rem}

%s3.2 #&#
\subsection{Asymptotic theory for partially linear single-index models}\label{sec3.2}

Denote $\vartheta_0 =(\beta_0^\top, \theta_0^\top)^\top$ and $\vartheta
=(\beta^\top, \theta^\top)^\top$ as a generic parameter for simplicity.
Let $\mathfrak{S}_n(\vartheta)$ and $\mathfrak{J}_n(\vartheta)$ be the
respective score and Hessian functions of $L_n(\vartheta)$ in the
minimisation problem \eqref{b9}. Let $\widehat{\vartheta}_n$ be the
estimator of $\vartheta_0$. We then have the expansion:
%
%e3.15 #&#
\begin{equation}
\label{c2a} 0=\mathfrak{S}_n(\widehat{\vartheta}_n)=
\mathfrak{S}_n(\vartheta_0 )+ \mathfrak{J}_n(
\vartheta_n) (\widehat{\vartheta}_n -\vartheta_0),
\end{equation}
where $\mathfrak{J}_n(\vartheta_n)$ is the Hessian matrix with the rows
evaluated at a point $\vartheta_n$ between $\widehat{\vartheta}_n$ and
$\vartheta_0$.

We also need to rotate our index vectors in model \eqref{a1}, namely,
reparametrerising the model, in order to derive the asymptotics. Using
the orthogonal matrix $Q=(\theta_0, Q_2)$ again, we can rewrite the
model as
%
%e3.16 #&#
\begin{eqnarray}
\label{c2a1} y_t=\beta_0^\top
QQ^\top x_t+g\bigl(\theta_0^\top
QQ^\top x_t\bigr)+e_t=\lambda
_{0}^\top z_{t}+ g\bigl(\alpha_{0}^\top
z_{t}\bigr)+e_t,
\end{eqnarray}
where $\lambda_0=Q^\top\beta_0=(\lambda_{10},\lambda_{20}^\top)^\top$
with $\lambda_{10}=\theta_0^\top\beta_0$ a scalar, $\lambda
_{20}=Q_2^\top\beta_0$ a $(d-1)$-dimensional vector, $\alpha
_{0}=Q^\top\theta_0$, $z_{t}=Q^\top x_t$ are defined the same as
before. Let $\lambda=Q^\top\beta$ and $\alpha=Q^\top\theta$ for the
generic vector rotation. Also, group them by $\mu_0= (\lambda_0^\top,
\alpha_0^\top)^\top$ and $\mu=(\lambda^\top, \alpha^\top)^\top$.

Let $L_n(\mu)$ be the counterpart of $L_n(\beta,\theta)$ after
reparameterisation. If $\widehat{\mu}_n$, the minimiser of $L_n(\mu)$,
is the estimator of $\mu_0$, then $\widehat{\mu}_n =\operatorname{diag} (Q^\top
, Q^\top) \widehat{\vartheta}_n$. Moreover, the score function
$\mathfrak{S}_n(\mu)$ and the Hessian $\mathfrak{J}_n(\mu)$ for the
parameter $\mu$ can be obtained from those for $\vartheta$. Namely,
$\mathfrak{S}_n (\mu)= \operatorname{diag} (Q^\top, Q^\top) \mathfrak
{S}_n(\vartheta)$ and $\mathfrak{J}_n (\mu)=\operatorname{diag}(Q^\top, Q^\top)
\mathfrak{J}_n (\vartheta) \operatorname{diag}(Q, Q)$. Premultiplying equation
\eqref{c2a} by $\operatorname{diag}(Q^\top, Q^\top)$, we have
%
%e3.17 #&#
\begin{equation}
\label{c2b} 0=\mathfrak{S}_n(\widehat{\mu}_n)=
\mathfrak{S}_n(\mu_0)+\mathfrak {J}_n(
\mu_n) (\widehat{\mu}_n -\mu_0),
\end{equation}
from which the following theorem is derived.

%th3.4 #&#
\begin{theo}\label{th4}
Under Assumptions \ref{assa} and \ref{assb}, $\widehat{\mu}_n\to_P\mu_0$. Moreover, as
$n\rightarrow\infty$
%
%e3.18 #&#
%e3.19 #&#
\begin{eqnarray}
n(\widehat{\lambda}_n-\lambda_0)&\to_D&
Q^\top \biggl( \int_0^1V(r)
V^\top (r)\,dr \biggr)^{-1} \int_0^1V(r)\,dU(r),
\label{th4a}
\\
D_n(\widehat{\alpha}_n-\alpha_0)
&\to_D& R^{-1/2}W(1), \label{th4b}
\end{eqnarray}
where $(U(r), V(r))$ is given in Assumption \ref{assa}, $D_n$, $R$ and $W$ are
the same as in Theorem~\ref{th31}.
\end{theo}

Theorem~\ref{th4} shows that for the partially linear single-index
model, the estimators of the parameters in the linear part have the
same rates of convergence as those in the linear model, while the
estimator of the index vector retains the dual rates in the system of
$Q$. Hence, there is a trio of rates of convergence accommodated in the
partially linear single-index model case. From Theorem~\ref{th4}, we
may derive asymptotic distributions for both $\widehat{\beta}_n$ and
$\widehat{\theta}_n$.

%th3.5 #&#
\begin{theo}\label{th5}
Under Assumptions \ref{assa} and \ref{assb}, for ($\widehat{\beta}_n$, $\widehat{\theta
}_n$) given by \eqref{b9} we have, as $n\rightarrow\infty$,
%
%e3.20 #&#
%e3.21 #&#
\begin{eqnarray}
n(\widehat{\beta}_n-\beta_0)&\to_D& \biggl(
\int_0^1V(r)V^\top(r)\,dr
\biggr)^{-1} \int_0^1V(r)\,dU(r),\label{th5a}
\\
n^{1/4}(\widehat{\theta}_n-\theta_0)
&\to_D&\mathbf{ MN}\bigl(0, \rho_{11}\theta_0
\theta_0^\top\bigr). \label{th5b}
\end{eqnarray}
\end{theo}

Furthermore, using \eqref{th4b}, for $\widehat{\theta}_n$, the results
of Theorems \ref{th32}--\ref{th33} and Corollaries \ref{cor1}--\ref
{cor2} with $\widehat{\theta}_{n, \mathrm{emp}}$ and $\widehat{g}_n(u)$
defined in the same fashion remain true.

%th3.6 #&#
\begin{theo}\label{th6}
Under Assumptions \ref{assa} and \ref{assb}, the results of
Theorems \ref{th32}--\ref
{th33} and Corollaries \ref{cor1}--\ref{cor2} also remain true for
$\widehat{\theta}_n$ and $\widehat{g}_n(u)$ in model \eqref{a1}.
\end{theo}

The proof of the main results in this section is given in Appendix \ref{appb}
below, except that Theorem~\ref{th31} and Corollaries \ref{cor1}--\ref{cor2} are shown in
Appendix D in the supplementary
material [\citet{DGT2015}].

%s4 #&#
\section{Extension to the general $H$-regular class}\label{sec4}

For integrated time series, the rate of convergence of the unknown
parameters involved in a regression function heavily depends on the
functional form of the regression function under consideration. The
literature focuses on two classes of functions, that is, the so-called
$I$-regular class and $H$-regular class. Integrable functions belong to
the $I$-regular class, while functions like power functions and
polynomial functions are $H$-regular. For more detail, we refer to
\citeauthor{phillips1999} (\citeyear{phillips1999,phillips2001}).

Observe that the partially linear single-index model is a combination
of the two classes. As already stated, this can be seen as a smooth
transition model whose behaviour in the finite domain is a linear model
perturbed by the $g$-function component, the influence of which is
reduced for a large $\|x_t\|$. Such models have broad applications.
Nonetheless, if the linear part may be relaxed to a nonlinear form,
the model will be more flexible and more applicable. To do so, we
combine a general $H$-regular function with a general $I$-regular
function to introduce a partially nonlinear single-index model of the form:
%
%e4.1 #&#
\begin{equation}
\label{ex1} y_t=f\bigl(\beta_0^\top
x_t\bigr)+g\bigl(\theta_0^\top
x_t\bigr)+e_t,\qquad  t=1,\ldots, n,
\end{equation}
where $f(\cdot)$ is parametrically known and $H$-regular, $g(\cdot)$ is
nonparametrically unknown and integrable, and $\beta_0$ and $\theta_0$
are the unknown parameters.

In some circumstances, one may have some idea on the trend in $y_t$
generated by an index variable $\beta^\top x_t$, for example, linear or
quadratic. Thus, model \eqref{ex1} should give an accurate description
for such a relation. Certainly, the partially linear single-index model
is a special case of model \eqref{ex1} with $f(u)=u$. The objective of
this section is then to estimate $(\beta_0, \theta_0)$ and $g(\cdot)$.

Before we propose our estimation method, we give a definition for the
$H$-regular class.

%de4.1 #&#
\begin{defn}\label{def}
We say that the function $f(x)$ is asymptotically homogeneous, or
$H$-regular, if for all $\eta>0$
%
%e4.2 #&#
\begin{equation}
\label{def1} f(\eta x)=\upsilon(\eta)F(x)+\xi(\eta; x), \qquad\bigl|\xi(\eta; x)\bigr|\le a(
\eta) P(x),
\end{equation}
where $F(x)$ and $P(x)$ are both locally integrable, and $\limsup_{\eta
\to\infty} \frac{a(\eta)}{ \upsilon(\eta)}=0$.
\end{defn}

If $f$ is $H$-regular with $\upsilon$ and $F$ satisfying \eqref{def1},
we call $\upsilon$ and $F$ the asymptotic order and the limit
homogeneous function of $f$, respectively. Note that any polynomial and
power function with positive power are $H$-regular. Note also that this
definition is not the exact one in the reference above, since in this
section $f(\cdot)$ is required to be differentiable.

\textit{Estimation procedure}: The estimation procedure follows similarly
from that for model \eqref{a1}. Using expansion \eqref{b3}, for each
$t$ rewrite \eqref{ex1} as $y_t-f(\beta^\top_0 x_t)=Z_k(\theta^\top_0
x_t)^\top c+ \gamma_k (\theta^\top_0 x_t)+e_t$, where $Z_k(\cdot)$, $c$
and $\gamma_k(\cdot)$ are defined as before. Let $\tilde{Y}= (y_1
-f(\beta^\top_0 x_1),\ldots, y_n-f(\beta^\top_0 x_n))^\top$, and $Z$,
$\gamma$ and $e$ remain the same as before. We then have the matrix
form equation, $\tilde{Y} =Zc+ \gamma+ e$. Then the OLS gives
$\widetilde{c}= \widetilde{c} (\beta_0,\theta_0) =(Z^\top Z)^{-1}Z^\top
\tilde{Y}$.

Define, $L_n(\beta, \theta)=\frac{1}{2}\sum_{t=1}^n [y_t-f(\beta^\top
x_t)-Z_k^\top(\theta^\top x_t)\widetilde{c} (\beta,\theta)]^2$. The
estimator of ($\beta_0$, $\theta_0$) is given by
%
%e4.3 #&#
\begin{equation}
\label{ex4} %
\pmatrix{\widehat{\beta}_n \vspace*{2pt}
\cr
\widehat{\theta}_n } %
=\mathop{\operatorname{argmin}}_{ \theta\in\Theta, \beta}
L_n(\beta, \theta).
\end{equation}
Similarly, a plug-in estimator $\widehat{g}_n(u)\equiv Z_k(u)^\top
\widehat{c}$ for any real $u\in\mathbb{R}$, where $\widehat
{c}=\widetilde{c}(\widehat{\beta}_n, \widehat{\theta}_n)$ is obtained
once $(\widehat{\beta}_n, \widehat{\theta}_n)$ is available.

\textit{Asymptotic theory}: The same notation as in Section~\ref{sec3.2} is used
for the minimisation problem \eqref{ex4}. Also, in order to derive the
corresponding asymptotic theory, we need to rotate vectors. Using the
orthogonal matrix $Q=(\theta_0, Q_2)$ again,
\begin{eqnarray*}
y_t=f\bigl(\beta_0^\top QQ^\top
x_t\bigr)+g\bigl(\theta_0^\top
QQ^\top x_t\bigr)+e_t=f\bigl(
\lambda_0^\top z_t\bigr)+ g\bigl(
\alpha_{0}^\top z_{t}\bigr)+e_t,
\end{eqnarray*}
where the notation used is the same as in \eqref{c2a1}, $\lambda
=Q^\top\beta$ and $\alpha=Q^\top\theta$ for generic vector rotation.
We also define $\mu_0= (\lambda_0^\top, \alpha_0^\top)^\top$ and $\mu
=(\lambda^\top, \alpha^\top)^\top$.

It is still true that if $\widehat{\mu}_n$ is the estimator of $\mu_0$
given by the minimiser of $L_n(\mu)$, then $\widehat{\mu}_n =\operatorname
{diag} (Q^\top, Q^\top) \widehat{\vartheta}_n$. Moreover, $\mathfrak
{S}_n (\mu)=\operatorname{diag}(Q^\top, Q^\top) \mathfrak{S}_n(\vartheta)$ and
$\mathfrak{J}_n (\mu)=\operatorname{diag}(Q^\top, Q^\top) \mathfrak
{J}_n(\vartheta) \operatorname{diag}(Q, Q)$.

%th4.1 #&#
\begin{theo}\label{thex}
Suppose that \textup{(i)} $f$ is $H$-regular with asymptotic order $\upsilon$
and limit homogeneous function $F$; \textup{(ii)} $\dot{f}$ and $\ddot{f}$
are
$H$-regular with asymptotic order $\dot{\upsilon}$ and limit
homogeneous function $\dot{F}$, and asymptotic order $\ddot{\upsilon}$
and limit homogeneous function $\ddot{F}$, respectively; \textup{(iii)} $|F(\beta
_0^\top x)-F(\beta^\top x)|$ is not a zero function on $\|x\| <\delta$
for some $\delta>0$ and if $\beta\neq\beta_0$; \textup{(iv)} $\upsilon(\sqrt
{n})^{-1} \sqrt{k}^3\to0$, where $k$ is the truncation parameter
satisfying Assumption \ref{assb}; (v)\break $|\dot{\upsilon}^{-2} (u)\ddot{\upsilon
}(u)\upsilon(u)|$ is bounded in $u\ge M_0$ for some $M_0>0$.

Under Assumptions \ref{assa} and \ref{assb}, we have $\widehat{\mu}_n\to_P\mu_0$.
Moreover, as $n\rightarrow\infty$
%
%e4.4 #&#
%e4.5 #&#
\begin{eqnarray}
n\dot{\upsilon}(\sqrt{n}) (\widehat{\lambda}_n-
\lambda_0)&\to_D& Q^\top \biggl( \int
_0^1\bigl[\dot{F}\bigl(\beta_0^\top
V(r)\bigr)\bigr]^2 V(r) V^\top(r)\,dr \biggr)^{-1}
\nonumber
\\[-8pt]
\\[-8pt]
\nonumber
&&{}\times\int_0^1\dot{F}\bigl(
\beta_0^\top V(r)\bigr)V(r)\,dU(r), %
\label{thexa}
\\
D_n(\widehat{\alpha}_n-\alpha_0)
&\to_D& R^{-1/2}W(1), \label{thexb}
\end{eqnarray}
where $(U(r), V(r))$, $D_n$, $R$ and $W$ are the same as in Theorem~\ref{th31}.
\end{theo}

%re4.1 #&#
\begin{rem}
It is reasonable to require that the derivatives of $f$ are
$H$-regular if $f$ is $H$-regular, as stated in conditions (i) and
(ii). Condition (iii) is simply an identification condition, while (iv)
and (v) are technical requirements that can be fulfilled easily by many
usual $H$-regular functions. Similar conditions for parameter
estimation in regression models involving $I(1)$ processes can be found
in \citet{phillips2001}. Particularly, $f(x)=x$ is a special case such
that conditions (i)--(v) are trivially satisfied.
\end{rem}

Similar to Theorems \ref{th5} and \ref{th6}, we derive some corresponding limit
distributions for $\widehat{\beta}_n$ and $\widehat{\theta}_n$ as well
as a plug-in estimate $\widehat{g}_n(u)$ below.

%th4.2 #&#
\begin{theo}\label{2thex}
Under the same conditions as Theorem~\ref{thex}, we have as
$n\rightarrow\infty$
%
%e4.6 #&#
%e4.7 #&#
\begin{eqnarray}
n\dot{\upsilon}(\sqrt{n}) (\widehat{\beta}_n-
\beta_0)&\to_D& \biggl( \int_0^1
\bigl[\dot{F}\bigl(\beta_0^\top V(r)\bigr)
\bigr]^2 V(r) V^\top(r)\,dr \biggr)^{-1}
\nonumber
\\[-8pt]
\\[-8pt]
\nonumber
&&{}\times\int_0^1\dot{F}\bigl(
\beta_0^\top V(r)\bigr)V(r)\,dU(r), %
\label{2thexa}
\\
n^{1/4}(\widehat{\theta}_n-\theta_0)
&\to_D& \mathbf{ MN}\bigl(0, \rho_{11}\theta
\theta_0^\top\bigr), \label{2thexb}
\end{eqnarray}
where the same notation is used as in Theorem~\ref{th5}.

Also, a plug-in estimate of the form: $\widehat{g}_n(u)=Z_k^\top
(u)\widehat{c}$ has the asymptotic normality as in Theorem~\ref{th33}. The
results in Theorem~\ref{th32} and Corollaries \ref{cor1}--\ref{cor2} remain valid.
\end{theo}

The proofs of Theorems \ref{thex}--\ref{2thex} are given in Appendix E of the supplementary
material
[\citet{DGT2015}].

%s5 #&#
\section{Simulation experiments}\label{sec5}

This section studies the finite-sample performance of the proposed
estimates. Let $d=2$ and $x_t$ be generated by
%
%e5.1 #&#
\begin{equation}
\label{4a} x_t=x_{t-1}+v_t \qquad\mbox{with }
v_t=r_0v_{t-1}+\varepsilon_t,
\end{equation}
for $t=1,\ldots, n$, where $r_0=0.1$, $\varepsilon_t\sim iiN(0, \sigma
^2 I_2)$, $x_0=0$ surely. Let sample size $n=400, 600$ and $1000$. The
number of Monte Carlo replications is $M=2000$. The truncation
parameter is $k=[a\cdot n^{\kappa}]$ with $\kappa=\frac{5}{44}$ and
$a=3.65$, satisfying the conditions in Assumption \ref{assb}. We shall then use
two examples.

%ex1 #&#
\begin{example}\label{ex5.1}
Consider a single-index model $y_t=
g(\theta_0^\top x_t)+ e_t$, $e_t\sim N(0,1)$, $t=1,\ldots, n$. There
are two parts in the simulation, according as $\theta_0^\top=
(0.6,-0.8)$ and $\theta_0^\top=(1,0)$ that both satisfy $\|\theta_0\|=1$.

We calculate the bias and standard deviation for $\widehat{\theta}_n=
(\widehat{\theta}_{1n}, \widehat{\theta}_{2n})^\top$:
%
%e5.2 #&#
\begin{equation}
\label{4b} \operatorname{Bias}=\bar{\widehat{\theta}}_n -
\theta_0,\qquad  \operatorname{S.d.}= \Biggl( \frac{1}{M} \sum
_{\ell=1}^M (\widehat{\theta}_{n\ell}
-\bar{\widehat {\theta}}_n)^{\otimes2} \Biggr)^{\otimes1/2},
\end{equation}
where $\otimes$ denotes an element-wise operation, and $\bar{\widehat
{\theta}}_n= \frac{1}{M} \sum_{\ell=1}^M \widehat{\theta}_{n\ell}$, in
which $\widehat{\theta}_{n\ell}$ stands for the $\ell$th replication of
the estimate.

In order to evaluate the asymptotic theory given in Theorem~\ref{th32},
we also calculate the bias and the standard deviation of $\widehat
{\theta}_{n,\mathrm{emp}}=\widehat{\theta}_n/\| \widehat{\theta}_n \|$ in
the same way as in~(\ref{4b}).

\textit{Part} I. Set $\theta_0^\top= (0.6, 0.8)$, $\sigma=0.6$ and $g(u)=
(1+u^2) e^{-u^2}$. We use the proposed procedure in Section~\ref{sec2.1} to
estimate $\theta_0$.

%t1 #&#
\begin{table}[b]
\caption{Bias and standard deviation for single-index model}
\label{sid1}
\begin{tabular*}{\textwidth}{@{\extracolsep{\fill}}lcd{2.4}d{2.4}d{2.4}cccc@{}}
\hline
& & \multicolumn{3}{c}{\textbf{Bias}} & & \multicolumn{3}{c@{}}{\textbf{S.d.}} \\[-6pt]
& & \multicolumn{3}{c}{\hrulefill} & & \multicolumn{3}{c@{}}{\hrulefill} \\
\multicolumn{1}{@{}l}{$\bolds{n}$}& & \multicolumn{1}{c}{\textbf{400}} & \multicolumn{1}{c}{\textbf{600}} & \multicolumn{1}{c}{\textbf{1000}} & & \multicolumn{1}{c}{\textbf{400}}
&\multicolumn{1}{c}{\textbf{600}} & \multicolumn{1}{c@{}}{\textbf{1000}} \\
\hline
$\widehat{\theta}_{1n}$ & &-0.0647 & -0.0519 &-0.0388 & &0.2678 &
0.2507 & 0.2042\\
$\widehat{\theta}_{2n}$ & &-0.0832 & -0.0684 &-0.0453 & &0.3461 &
0.3285 & 0.2586\\
$\widehat{\theta}_{n,\mathrm{emp}}^1$ & &0.0043 & 0.0024 & -0.0016 &
&0.1005 &0.0820 & 0.0679 \\
$\widehat{\theta}_{n,\mathrm{emp}}^2$ & &0.0063 & 0.0066 &0.0050 &
&0.0717 &0.0659 & 0.0515 \\
\hline
\end{tabular*}
\end{table}

As can be seen from Table~\ref{sid1}, both the biases and the standard
deviations for $\widehat{\theta}_n$ decrease as the sample size
increases, and $\widehat{\theta}_{1n}$ and $\widehat{\theta}_{2n}$ have
similar performance. Moreover, the biases and standard deviations of
$\widehat{\theta}_{n, \mathrm{emp}}$ indicate that $\widehat{\theta}_{n,
\mathrm{emp}}$ has a rate of convergence faster than that of $\widehat
{\theta}_n$, as shown in Theorem~\ref{th32}.

\textit{Part} II. Put $\theta_0^\top= (1,0)$, $\sigma=0.6$ and
$g(u)=(1+u^2)\exp(-u^2)$. As pointed out before, the rotation of the
parameters is only for the derivation of the asymptotic theory. To
evaluate the asymptotic theory given in Theorem~\ref{th31}, we directly take
$\theta_0^{\top}=\alpha_0^{\top} =(1,0)$ so that $\widehat{\alpha}_n
=\widehat{\theta}_n$ in this experiment.

%t2 #&#
\begin{table}
\caption{Bias and standard deviation for single-index model}
\label{sid2}
\begin{tabular*}{\textwidth}{@{\extracolsep{\fill}}lcd{2.4}d{2.4}d{2.4}cccc@{}}
\hline
& & \multicolumn{3}{c}{\textbf{Bias}} & & \multicolumn{3}{c@{}}{\textbf{S.d.}} \\[-6pt]
& & \multicolumn{3}{c}{\hrulefill} & & \multicolumn{3}{c@{}}{\hrulefill} \\
\multicolumn{1}{@{}l}{$\bolds{n}$}& & \multicolumn{1}{c}{\textbf{400}} & \multicolumn{1}{c}{\textbf{600}} & \multicolumn{1}{c}{\textbf{1000}} & & \multicolumn{1}{c}{\textbf{400}}
&\multicolumn{1}{c}{\textbf{600}} & \multicolumn{1}{c@{}}{\textbf{1000}} \\
\hline
$\widehat{\alpha}_{1n}$ & &0.0866 & 0.0768 &0.0340 & &0.3803 & 0.3748 &
0.3338\\
$\widehat{\alpha}_{2n}$ & &0.0013 & -0.0008 &-0.0006 & &0.1388 & 0.1186
& 0.0898\\
$\widehat{\alpha}_{n,\mathrm{unit}}^1$ & &-0.0073 & -0.0061 & -0.0031 &
&0.0246 &0.0237 & 0.0128 \\
$\widehat{\alpha}_{n,\mathrm{unit}}^2$ & &0.0011 & -0.0018 &-0.0003 &
&0.1186 &0.1080 & 0.0779 \\
\hline
\end{tabular*}
\end{table}

As can be seen from Table~\ref{sid2}, both the biases and the standard
deviations of $\widehat{\alpha}_{1n} $ and $\widehat{\alpha}_{2n}$
decrease as the sample size increases. Particularly, the decrease for
$\widehat{\alpha}_{2n}$ is much faster than that for $\widehat{\alpha
}_{1n}$. This verifies the type of rates of convergence given in
Theorem~\ref{th31} that $\widehat{\alpha}_{2n}- \alpha
_{20}=O_P(n^{-3/4})$, while $\widehat{\alpha}_{1n}- \alpha_{10}= O_P(n^{-1/4})$.

Nonetheless, shown by the standard deviations, $\widehat{\alpha}_{n,
\mathrm{unit}}^1$ converges significantly faster than $\widehat{\alpha
}_{n, \mathrm{unit}}^2$. This is also implied by Corollary~\ref{cor1}
that $\widehat{\alpha}_{n, \mathrm{unit}}^1-\alpha_{10}=O_P(n^{-3/2})$
and $\widehat{\alpha}_{n, \mathrm{unit}}^2 -\alpha_{20}=O_P(n^{-3/4})$.
Note also that the biases of $\widehat{\alpha}_{n, \mathrm{unit}}^1$ are
always negative (by definition, $\widehat{\alpha}_{n, \mathrm{unit}}^1\le
\alpha_{10}=1$) for each Monte Carlo experiment. As a result, the
biases of $\widehat{\alpha}_{n, \mathrm{unit}}^1$ approach zero
relatively slower than those of $\widehat{\alpha}_{n, \mathrm{unit}}^2$.
In addition, $\widehat{\alpha}_{n, \mathrm{unit}}^2$ and $\widehat{\alpha
}_{2n}$ perform quite similarly since they have the same rate of convergence.
\end{example}

%ex2 #&#
\begin{example}\label{ex5.2}
In this example, a partially linear
single-index model of the form: $y_t= \beta_0^\top x_t+g(\theta_0^\top
x_t)+ e_t$, $e_t\sim N(0,1)$, $t=1,\ldots, n$, is examined with
$g(u)=(1+u^2) \exp(-u^2)$, $\beta_0^\top=(0.3,0.5)$, $\theta_0^\top
=(0.6, -0.8)$ and $\sigma=0.8$ involved in $\varepsilon_t \sim i i N(0,
\sigma^2  I_2)$.

%t3 #&#
\begin{table}
\caption{Bias and standard deviation for partially linear single-index model}
\label{plsid}
\begin{tabular*}{\textwidth}{@{\extracolsep{\fill}}lcd{2.4}d{2.4}d{2.4}cccc@{}}
\hline
& & \multicolumn{3}{c}{\textbf{Bias}} & & \multicolumn{3}{c@{}}{\textbf{S.d.}} \\[-6pt]
& & \multicolumn{3}{c}{\hrulefill} & & \multicolumn{3}{c@{}}{\hrulefill} \\
\multicolumn{1}{@{}l}{$\bolds{n}$}& & \multicolumn{1}{c}{\textbf{400}} & \multicolumn{1}{c}{\textbf{600}} & \multicolumn{1}{c}{\textbf{1000}} & & \multicolumn{1}{c}{\textbf{400}}
&\multicolumn{1}{c}{\textbf{600}} & \multicolumn{1}{c@{}}{\textbf{1000}} \\
\hline
$\widehat{\theta}_{1n}$ & & -0.0495 &-0.0470 &-0.0324 & & 0.2652
&0.2494 &0.1991 \\
$\widehat{\theta}_{2n}$ & &0.0676 &0.0645 &0.0435 & &0.3433 & 0.3340
&0.2572 \\
$\widehat{\theta}_{n,\mathrm{emp}}^1$ & &0.0038 & 0.0031 &0.0023 &
&0.0934 &0.0798 & 0.0597\\
$\widehat{\theta}_{n,\mathrm{emp}}^2$ & &-0.0062 & -0.0041 &-0.0019 &
&0.0761 &0.0621 &0.0475 \\
$\widehat{\beta}_{1n}$ & &-0.0010 & -0.0002 &0.0001 & & 0.0106 &0.0068
&0.0038 \\
$\widehat{\beta}_{2n}$ & &-0.0007 & 0.0001 &0.0001 & &0.0118 & 0.0067
&0.0037 \\
\hline
\end{tabular*}
\end{table}

Formulae in \eqref{4b} are used for $\widehat{\theta}_n$, $\widehat
{\theta}_{n, \mathrm{emp}}$ and $\widehat{\beta}_n$. All simulation
results with sample size $n=400, 600, 1000$ and $\sigma=0.8$ are
reported in Table~\ref{plsid}. As can be seen, both the biases and the
standard deviations decrease as the sample size increases. Moreover,
the rate of $\widehat{\theta}_{n, \mathrm{emp}}$ approaching the true
value looks faster than that of $\widehat{\theta}_n$. This is supported
by Theorems \ref{th5} and \ref{th6} that $\widehat{\theta}_n- \theta
_0=O_P(n^{-1/4})$ and $\widehat{\theta}_{n, \mathrm{emp}} -\theta
_0=O_P(n^{-3/4})$.

%f1 #&#
\begin{figure}[b]

\includegraphics{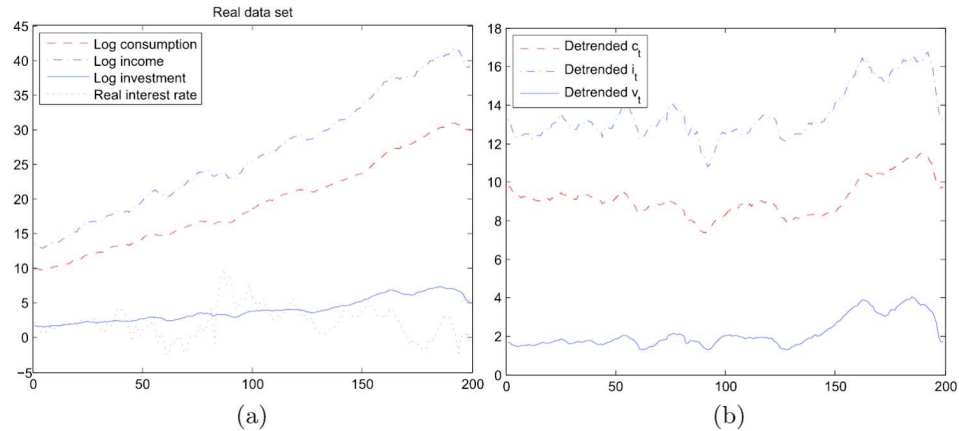}

\caption{The real data and the de-trended data.
\textup{(a)} The real data. \textup{(b)} The detrended data.}
\label{5f1}
\end{figure}

Meanwhile, since $\widehat{\beta}_n= (\widehat{\beta}_{1n}, \widehat
{\beta}_{2n})^\top$ possesses the fastest rate of convergence of
$n^{-1}$ by Theorem~\ref{th5}, both the biases and the standard
deviations of $\widehat{\beta}_n$ support the large sample behaviour.
Therefore, the asymptotic theory established in Section~\ref{sec3} has been
evaluated in these examples.
\end{example}

%s6 #&#
\section{Empirical study}\label{sec6}

We propose to use a partially linear single-index model to fit an
empirical data set before we make some comparisons with some candidate models.

\textit{The data}. The aggregate US data on consumption,
income, investment and interest rate are obtained from \emph{Federal
Reserve Economic Data (FRED)}. We consider a quarterly data set over
1960:1--2009:3 with 199 observations. Let $r_t$ stand for the real
interest rate, and $c_t=\log(C_t), i_t=\log(I_t)$ and $v_t= \log
(V_t)$, where $C_t, I_t$ and $V_t$ are the consumption expenditures,
disposable incomes and investments, respectively, for $t=1,\ldots,
199$. The data of $c_t,i_t,v_t$ and $r_t$ are plotted in (a) of Figure~\ref{5f1}. It can be seen that all of them have trending components
except $r_t$. To meet the theoretical assumptions, we de-trend the data
for $c_t, i_t$ and $v_t$. More precisely, suppose that $c_t=\mu
_1+c_{t-1}+u_{1t}$, $i_t=\mu_2+ i_{t-1}+u_{2t}$ and $v_t=\mu_3 +v_{t-1}
+u_{3t}$ for $t=2,\ldots,199$, where $u_{it}$, $i=1,2,3$, are error
terms. Then $\mu_i$ are estimated as: $\widehat{\mu}_1=\frac{1}{198}
\sum_{i=2}^{199} (c_t- c_{t-1})=0.1022$, $\widehat{\mu}_2= \frac
{1}{198} \sum_{i=2}^{199}(i_t-i_{t-1}) =0.1302$, $\widehat{\mu}_3= \frac
{1}{198} \sum_{i=2}^{199} (v_t-v_{t-1}) =0.0181$.

Define for each $t$, $\tilde{c}_t=c_t-\widehat{\mu}_1t$, $\tilde
{i}_t=i_t -\widehat{\mu}_2t$ and $\tilde{v}_t=v_t-\widehat{\mu}_3t$.
They are the de-trended versions being plotted in (b) of Figure~\ref{5f1}, correspondingly.

%\begin{figure}[htp]
%
%\centering
%\subfloat[]{\includegraphics[width=2.4in]{realciii.pdf}}
%%\label{realdata}
%\subfloat[]{\includegraphics[width=2.4in]{ddata.pdf}}
%%\label{ddata}
%\end{figure}

\begin{figure}[b]

\includegraphics{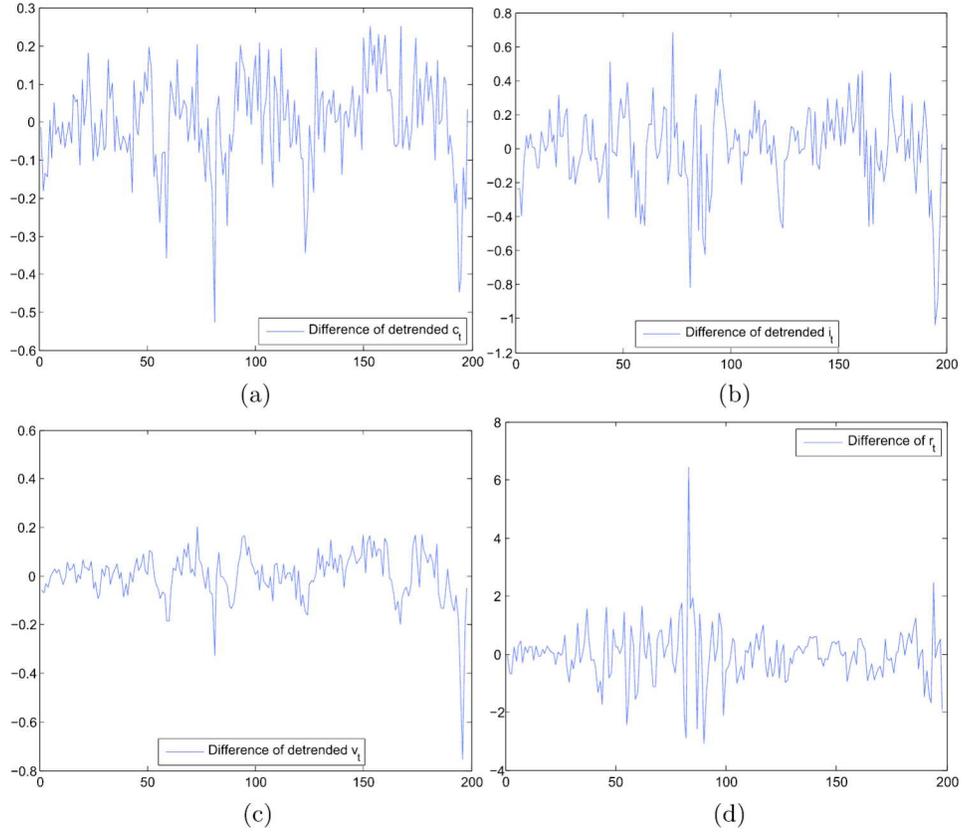}

\caption{Difference of dataset. \textup{(a)} Difference of detrended consumption.
\textup{(b)} Difference of detrended income.
\textup{(c)} Difference of detrended investment.
\textup{(d)} Difference of interest rate.}
\label{5f2}
\end{figure}

An ADF test is applied to each of $\tilde{c}_t$, $\tilde{i}_t$ and
$\tilde{v}_t$, respectively. The ADF test fails to reject the null of
possessing a unit root with $p$-values 0.7595, 0.6293 and 0.7637,
respectively. In addition, it is known that $r_t$ is an $I(1)$ process
[\citet{jiti2009a}]. To visualise the $I(1)$ processes, the plots of the
differences are given in Figure~\ref{5f2}.

%f2 #&#

%
%\begin{figure}[ht]
%\centering
%
%\subfloat[]{\includegraphics
%[width=2.4in]{dc.pdf}}
%\subfloat[]{\includegraphics
%[width=2.4in]{di.pdf}}\\
%\subfloat[]{\includegraphics
%[width=2.4in]{dv.pdf}}
%\subfloat[]{\includegraphics[width=2.4in]{dr.pdf}}
%\end{figure}

\textit{The model}.  A partially linear single-index model is
proposed to fit the data $\tilde{c}_t $, $\tilde{i}_t$, $\tilde{v}_t$
as well as $r_t$ in the following forms:
%
%e6.1 #&#
\begin{equation}
\label{5c} y_t=\beta_0^\top
x_t+g\bigl(\theta_{0}^\top x_t
\bigr)+e_{t},
\end{equation}
where $t=2,\ldots, 199$, $y_t=\tilde{c}_t$ and $x_t^\top=(x_{1t},
x_{2t}, x_{3t}, x_{4t}, x_{5t})$ in which $x_{1t}=\tilde{i}_{t-1}$,
$x_{2t}=\tilde{i}_{t}$, $x_{3t}= \tilde{v}_t$, $x_{4t}=\tilde{v}_{t-1}
$, $x_{5t} =r_t$, and $g(\cdot)$ is an unknown integrable function,
$e_{t}$ is the error term. Note that we only include the first lagged
information in the discussion, as they are more relevant than the other lags.

%
%f3 #&#
\begin{figure}[b]

\includegraphics{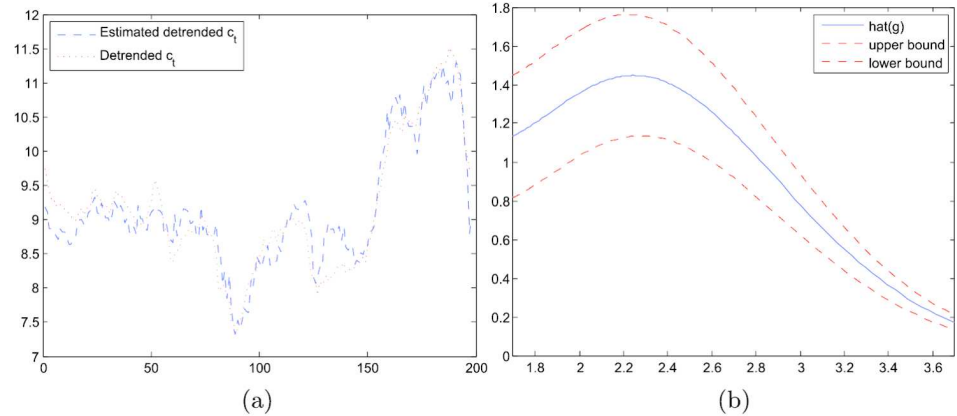}

\caption{Estimated data and estimated link function.
\textup{(a)} Model \protect\eqref{5c}.
\textup{(b)} Confidence interval curve.}
\label{5f3}
\end{figure}

\textit{Estimation}.  Before implementing our proposed
procedures to estimate model \eqref{5c}, one issue is to determine a
suitable truncation parameter $k$ so that the function $g(\cdot)$ can
be better approximated by the first $k$ terms in $\{\mathscr{H}_i(x)\}
$. Toward this end, we propose using the Generalised Cross-Validation
(GCV) method [see \citet{jiti2002}] as an initial step to select an
optimal value $k$. Note that while there is no theory for such
selection in the nonstationary time series case, the initial selection
method works numerically in this example. Let $\widehat{k}$ denote an
optimal value such that
%
%e6.2 #&#
\begin{equation}
\label{5e1} \widehat{k}=\mathop{\operatorname{argmin}}_{k\in K_n} \biggl(1-
\frac{k}{n} \biggr)^{-2} \widehat{\sigma}^2(k),
\end{equation}
where $\widehat{\sigma}^2(k) =\frac{1}{n} \sum_{t=1}^n  (y_t - \beta
^\top x_t- Z_k(\theta^\top x_t)^\top\widetilde{c}(\beta,\theta)
)^2$, $K_n=\{1,\ldots, 12\}$.

We have $\widehat{k}=5$ by GCV, $\widehat{\beta}_n =( -0.0479, 0.5701,
-1.1689, 1.8685, -0.1223)$ and $\widehat{ \theta}_{n} = ( 0.2110,
-0.3452, 0.0835, 2.6095, -0.2022)$. Meanwhile, we have $\widehat{c}=
\widetilde{c} (\widehat{\beta}_n, \widehat{\theta}_{n} ) = (-89.64,
112.54, -74.65, 28.94, -3.33)^\top$. This suggests
%
%e6.3 #&#
\begin{eqnarray}
\label{5g} %
\widehat{g}_{5}(u)&=&\bigl[-89.64
d_0^{-1} H_0(u)+112.54 d_1^{-1}H_1(u)-74.65
d_2^{-1}H_2(u)
\nonumber
\\[-8pt]
\\[-8pt]
\nonumber
&&{}+ 28.94 d_3^{-1}H_3(u) -3.33
d_4^{-1}H_4(u)\bigr]e^{-u^2/2},
\end{eqnarray}
where $H_i(u)$ are the Hermite polynomials, and $d_i=(\sqrt{\pi}
2^ii!)^{1/2}$ are the norm of $H_i(u)$ in $L^2(\mathbb{R}, e^{-u^2})$
for $i=0,1,\ldots, 4$.

In comparison, the de-trended log consumption $y_t=\tilde{c}_t$ is
plotted along with the estimated de-trended log consumption by the
partially linear single-index model $\widehat{y}_t= \widehat{\beta
}_n^\top x_t+ \widehat{g}_{5} ( \widehat{\theta}_{n}^\top x_t)$ in (a)
of Figure~\ref{5f3}.

%\begin{figure}[ht]
%
%\centering
%\subfloat[]{\includegraphics[width=2.4in]{plsgid.pdf}}
%%\subfloat[Model \eqref{5a}]{\includegraphics[width=2.5in]{linear.pdf}}\
%%\
%%\label{linear}
%%\subfloat[Model \eqref{5b}]{\includegraphics[width=2.5in]{sgid.pdf}}
%%\label{sgid}
%\subfloat[]{\includegraphics
%[width=2.4in]{gintval2.pdf}}
%\end{figure}

Note also by \eqref{5g} that the estimated link function $\widehat
{g}_{5}(u)$ is integrable on $\mathbb{R}$. According to the normality
in Theorem~\ref{th33}, we draw the confidence bands for $\widehat
{g}_{5}(u)$ at the significance level of 80\% in (b) of Figure~\ref{5f3}.

\textit{Comparison}. To check whether the estimated
relationship by the partially linear model is a suitable one, we shall
compare model \eqref{5c} with two natural competitors of the form
%
%e6.4 #&#
%e6.5 #&#
\begin{eqnarray}
y_t&=&h\bigl(\theta_{10}^\top x_t
\bigr)+e_{1t},\label{5a}
\\
y_t&=&\beta_{\mathrm{linear}}^\top x_t+e_{2t},
\label{5b}
\end{eqnarray}
where $h(\cdot)$ is integrable and unknown, and $e_{\ell t}$ are the
error terms for $\ell=1,2$.

To begin with, the linear model is estimated by OLS with $\widehat
{\beta}_{\mathrm{linear}}=( -0.0628,\break  0.7952, -1.2315, 0.9414,
-0.0644)^\top$. Moreover, GCV is applied for model \eqref{5a} with
$\widehat{\sigma}^2 (k) =\frac{1}{n} \sum_{t=1}^n  (y_t - Z_k(\theta
^\top x_t)^\top\widetilde{c}(\theta) )^2$, and we have $\widehat
{k}=3$. Then $\widehat{\theta}_{1n}=(-0.0014, 0.0152, -0.0229, 0.0176,
-0.0016)^\top$. Meanwhile, the estimate of $\widehat{c}= \widetilde
{c}(\widehat{\theta}_{1n})= (237.05, -61.92, 315.32)^\top$ implies
%
%e6.6 #&#
\begin{eqnarray}\qquad
\label{5f}&& \widehat{h}_{3}(u)
\nonumber
\\[-8pt]
\\[-8pt]
\nonumber
&&\qquad=\bigl[237.05 d_0^{-1}H_0(u)-61.92
d_1^{-1}H_1(u) +315.32d_2^{-1}H_2(u)
\bigr]e^{-u^2/2},
\end{eqnarray}
using the same notation as in \eqref{5g}.

To proceed further, we compare the so-called in-sample and
out-of-sample mean square errors among the three models.

(i) In-sample mean square error ($\mathrm{MSE}_{\mathrm{in}}$): As above, all
unknown parameters and functions in the three models \eqref{5c}, \eqref
{5a} and \eqref{5b} are estimated based on the whole observations
($x_t$, $y_t$), $t=2,\ldots, 199$. Once these have been done, we shall
have estimated $\widehat{y}^\ell_t$ with $\ell=1,2,3$ corresponding to
models \eqref{5c}, \eqref{5a}, \eqref{5b} for $t=2,\ldots, 199$,
\begin{eqnarray*}
\widehat{y}_t^1=\widehat{\beta}^\top_n
x_t+\widehat{g}_{5} \bigl(\widehat {
\theta}_{n}^\top x_t\bigr),\qquad
\widehat{y}_t^2=\widehat{h}_{3}\bigl(
\widehat{\theta}_{1n}^\top x_t\bigr)\quad \mbox{and}\quad
\widehat{y}_t^3=\widehat{\beta}_{\mathrm{linear}}^\top
x_t.
\end{eqnarray*}

Then the in-sample mean square errors are calculated, for $\ell=1,2,3$, by
%
%e6.7 #&#
\begin{equation}
\label{5d} \mathrm{MSE}_{\mathrm{in}}(\ell)=\frac{1}{198}\sum
_{t=2}^{199}\bigl(y_t-\widehat
{y}^\ell_t\bigr)^2.
\end{equation}

Meanwhile, to verify the choice of GCV, the $\mathrm{MSE}_{\mathrm{in}}$ for model
\eqref{5c} with $k=3,4,6,7$, respectively, and for model \eqref{5a}
with $k=1,2, 4,5$, respectively, are calculated as well.

(ii) Out-of-sample mean square error ($\mathrm{MSE}_{\mathrm{out}}$): Each time,
one part of observations is used to estimate all unknown parameters and
functions in the three models; then the next value of the dependent
variable is forecasted using the estimated models. More precisely,
letting $j=1, 2, \ldots, 10$, we use the observations $\{(y_t, x_t):
2\le t\le178+2j\}$ to estimate the unknown parameters and functions
[with fixed $\widehat{k}=5$ for model \eqref{5c} and $\widehat{k}=3$
for model \eqref{5a}] in the three models, then the next $y_{179+2j}$
is forecasted by the three models using the estimated parameters,
\begin{eqnarray*}
\widehat{y}_{179+2j}^1&=&\widehat{\beta}^\top_j
x_{179+2j}+ \widehat {g}_{5}^j\bigl( \widehat{
\theta}_{j}^\top x_{179+2j}\bigr),
\\
\widehat{y}_{179+2j}^2&=&\widehat{h}_{3}^j
\bigl(\widehat{\theta}_{1j}^\top x_{179+2j}\bigr)\quad
\mbox{and} \quad\widehat{y}_{179+2j}^3=\widehat{\beta
}_{j,\mathrm{linear}}^\top x_{179+2j}.
\end{eqnarray*}

The $\mathrm{MSE}_{\mathrm{out}}$ are evaluated, for $\ell=1,2,3$, by
%
%e6.8 #&#
\begin{equation}
\label{5e} \mathrm{MSE}_{\mathrm{out}}(\ell)=\frac{1}{10}\sum
_{j=1}^{10}\bigl(y_t-
\widehat{y}^\ell _ {179+2j}\bigr)^2.
\end{equation}

In addition, to assess the choice of GCV, the $\mathrm{MSE}_{\mathrm{out}}$ for
model \eqref{5c} with $k=3,4,6,7$, respectively, and for model \eqref
{5a} with $k=1,2, 4,5$, respectively, are computed as well. All $\mathrm
{MSE}_{\mathrm{in}}$ and $\mathrm{MSE}_{\mathrm{out}}$ are given
in Table~\ref{comparison}.
%
%t4 #&#
\begin{table}
\caption{The MSEs for models \protect\eqref{5c}, \protect\eqref{5a}, \protect\eqref{5b}}
\label{comparison}
\begin{tabular*}{\textwidth}{@{\extracolsep{\fill}}lcccccccc@{}}
\hline
$\bolds{k=}$ & &\textbf{3} &\textbf{4} &\textbf{5}& \textbf{6} & \textbf{7}& & \\
\hline
& & \multicolumn{5}{c@{}}{Partially linear single-index model \eqref{5c}}\\
%\\ [-6pt]
%& & \multicolumn{5}{c@{}}{\hrulefill}\\
$\mathrm{MSE}_{\mathrm{in}}$& &0.0968 & 0.1018 & 0.0946 &0.1460 &0.1559\\
$\mathrm{MSE}_{\mathrm{out}}$& &0.2418 & 0.1761 & 0.1232 & 0.2146 &0.1786\\ \hline
& & \multicolumn{5}{c}{Single-index model \eqref{5a}} & & \multicolumn{1}{c@{}}{Linear model
\eqref{5b}}\\
%& & \multicolumn{5}{c}{\hrulefill} & &
%\multicolumn{1}{c@{}}{\hrulefill}\\
%$\bolds{k=}$ & &\textbf{3} &\textbf{4} &\textbf{5}& \textbf{6} & \textbf{7}& & \\
%\hline
$\mathrm{MSE}_{\mathrm{in}}$ & &0.3011 & 0.1641& 0.1544 & 0.7709 &1.4076 & &
0.1666 \\
$\mathrm{MSE}_{\mathrm{out}}$ & &0.6962 & 0.2733 & 0.2607 & 1.9060 &3.0838 & &
0.2598\\
\hline
\end{tabular*}
\end{table}

In summary, among the three models, the partially linear single-index
model \eqref{5c} performs much better than the other two, in the sense
that both its $\mathrm{MSE}_{\mathrm{in}}$ and $\mathrm{MSE}_{\mathrm{out}}$ are the smallest
within the models. Particularly, model \eqref{5c} outperforms models
\eqref{5a} and \eqref{5b} over all choices of the truncation parameter
regardless of whether or not it is chosen by GCV method. This is
possibly because model \eqref{5c} is the combination of a linear trend
and a local adjustment by the link function such that it is more
flexible than the other two.

Note also that, with $\widehat{k}=5$ model \eqref{5c} has the best
performance. Therefore, model \eqref{5c} with $\widehat{k}=5$ is the
most favourable one to explain the empirical relationship between the
consumption and the income, investment and real interest rate for the
US data from the period of 1960 to 2009.

%s7 #&#
\section{Conclusions}\label{sec7}

The estimation procedures for both single-index and partially
single-index models in the presence of nonstationarity and
integrability have been proposed. New asymptotic properties for the
proposed estimates have been established. The rate of convergence of
the estimators of the index vector $\theta_0$ consists of two different
components in a new coordinate system for both the single-index and the
partially linear single-index models, while the estimator of the
coefficient vector $\beta_0$ has the super $n$-rate. The normality of
the plug-in estimate of the link function involved in each model has
been established. To satisfy the identification condition, the
normalisation of the estimator of $\theta_0$ in each case has been
proposed and interestingly it possesses a fast rate of convergence.
Motivated by more applicability, the partially linear single-index
model is extended by using a general $H$-regular function to replace
the linear function. New results have been obtained. Meanwhile, Monte
Carlo simulations have supported the key theoretical properties.
Furthermore, the empirical study has shown that the partially linear
single-index model outperforms both the linear and the single-index
models, and is the most suitable one for the aggregate US data on
consumption, income, investment and interest rate.

%sA #&#
\begin{appendix}\label{app}

%sA #&#
\section{Lemmas}\label{appa}

Three lemmas are given in this appendix while their proofs are shown in
Appendix C of the supplementary
material [\citet{DGT2015}].

%leA.1 #&#
\begin{lemma}\label{lemma5a}
The following assertions hold:
\begin{longlist}[(1)]
\item[(1)]$\frac{1}{\sqrt{t}}(x_{1t}, x_{2t}^\top)$ has a joint probability
density $\psi_t(x, w^\top)$; and given $\mathcal{F}_{s}$ (defined in
Assumption \ref{assa}), $\frac{1}{ \sqrt{t-s}} (x_{1t}-x_{1s}, x_{2t}^\top
-x_{2s}^\top)$ has a joint density $\psi_{ts}(x, w^\top)$ where
$t>s+1$. Meanwhile, these density functions are bounded uniformly in
$(x,w)$ as well as $t$ and $(t,s)$, respectively.

\item[(2)] For large $t$ and $t-s$, we have $\psi_t(x, w^\top)=\phi
(w)f_t(x)(1+o(1))$ and $\psi_{ts}(x, w^\top)= \phi(w)f_{ts}(x)
(1+o(1))$ where $\phi(w)$ is the density of an $(d-1)$-dimensional
normal distribution, $f_t(x)$ is the marginal density of $\frac{1}{\sqrt
{t}}x_{1t}$ and $f_{ts}(x)$ is the marginal density of $\frac{1}{\sqrt
{t-s}}(x_{1t}-x_{1s})$.
\end{longlist}
\end{lemma}

%leA.2 #&#
\begin{lemma}\label{lemma2}
\textup{(1)} Under Assumptions \ref{assa} and \ref{assb}, we have as $n\to\infty$, \\ $\llVert
\frac{1}{\sqrt{n}}Z^\top Z-L_1(1,0)I_{k} \rrVert  =o_P(1)$ in a richer
probability space.\vspace*{-6pt}
\begin{longlist}[(2)]
\item[(2)] Let $\widehat{Z}$ be the matrix $Z$ defined in Section~\ref{sec2} with
$\theta$ being replaced by $\widehat{\theta}_n$. Under Assumptions \ref{assa}
and \ref{assb}, we have $\frac{1}{\sqrt{n}}\llVert Z^\top Z- \widehat{Z}^{\top}
\widehat{Z}\rrVert  =o_P(1)$.
\end{longlist}
\end{lemma}

%leA.3 #&#
\begin{lemma}\label{cghat}
Under Assumptions \ref{assa} and \ref{assb}, we have $\|\widetilde{c}(\theta_0)-c\|^2
=O_P(1)\frac{k}{\sqrt{n}}$ as $n\to\infty$ where $\widetilde{c}(\theta
_0)$ is defined in Section~\ref{sec2.1}. %It is also true for $\widetilde{c}(
%\beta,\theta)$ defined in Section~2.2, viz., $\|\widetilde{c}(\beta,
%\theta) -c\|^2 =O_P(1) \frac{k}{\sqrt{n}}$ for $\{(\beta,\theta):\ \|
%\beta-\beta_0\|+\|\theta-\theta_0\|<\epsilon\}$ and small $\epsilon$.
\end{lemma}

%sA #&#
\section{Proofs of the main results}\label{appb}

The full proof of Theorem~\ref{th32} and the outlines of the proofs of
Theorems \ref{th33} and \ref{th4} are given below. In the meantime, all detailed
proofs for the theorems and corollaries in Section~\ref{sec3} and that in
Section~\ref{sec4} are given in Appendices D and E, respectively, of the
supplementary material [\citet{DGT2015}].

\begin{pf*}{Proof of Theorem~\ref{th32}}
Noting that $\sqrt[4]{n}D_n^{-1}\to\mathrm{diag}(1, \mathbf{0}_{d-1})$
as $n\to\infty$ where $\mathbf{0}_{d-1}$ is a zero matrix of
$(d-1)\times(d-1)$, by the continuous mapping theorem and Theorem~\ref{th31}
we have
%
%eA.1 #&#
\begin{eqnarray}
&&\sqrt[4]{n} (\widehat{\theta}_n-\theta_0)\nonumber\\
&&\qquad=
\sqrt[4]{n} \bigl(D_nQ^\top \bigr)^{-1}D_nQ^\top(
\widehat{\theta}_n -\theta_0)=Q\sqrt[4]{n}D_n^{-1}D_n
(\widehat{\alpha}_n -\alpha_0)
\\
&&\qquad\to_D Q \operatorname{diag}(1, \mathbf{0}_{d-1})R^{-1/2}W(1)=
\mathbf {MN}\bigl(0,\rho_{11} \theta_0 \theta_0^\top
\bigr). \nonumber%
\end{eqnarray}
In addition, it follows from $\widehat{\theta}_{n, \mathrm{emp}}=
Q\widehat{\alpha}_{n, \mathrm{unit}}$, $\theta_0=Q\alpha_0$ and Corollary~\ref{cor1} that
\begin{eqnarray*}
n^{3/4} (\widehat{\theta}_{n, \mathrm{emp}}-\theta_0 )
&=&Qn^{3/4} %
\pmatrix{\widehat{\alpha}_{n,\mathrm{unit}}^1-1
\vspace*{2pt}
\cr
\widehat{\alpha}_{n,\mathrm{unit}}^2 } %
=(
\theta_0 Q_2) %
\pmatrix{0 \vspace*{2pt}
\cr
n^{3/4}\widehat{\alpha}_{n,\mathrm{unit}}^2 }
+o_P(1)
\\
&=& Q_2n^{3/4} \widehat{\alpha}_{n,\mathrm{unit}}^2+o_P(1)
\to_D  \mathbf{ MN}\bigl(0,Q_2
\rho_{22}Q_2^\top\bigr).
\end{eqnarray*}
\upqed\end{pf*}

\begin{pf*}{Outline of the proof of Theorem~\ref{th33}}
The uniform consistency of $\widehat{g}_n(u)$ follows from Lemma~\ref
{cghat} directly. Indeed, for large $n$ and by the consistency of
$\widehat{\theta}_n$ and the continuity of $\widetilde{c}(\theta)$ in
$\theta$, we have $\|\widehat{c}-c\|^2=\|\widetilde{c}(\widehat{\theta
}_n)-c\|^2=O_P(1)\frac{k}{\sqrt{n}}$.
\begin{eqnarray*}
\sup_{u\in\mathbb{R}}\bigl|\widehat{g}_n(u)-g(u)\bigr|&\le&\sup
_{u\in\mathbb
{R}}\bigl|Z_k^\top(u)[\widehat{c}-c]\bigr|+\sup
_{u\in\mathbb{R}}\bigl|\gamma_k(u)\bigr| \le\sup_{u\in\mathbb{R}}
\bigl\|Z_k (u)\bigr\| \|\widehat{c} -c\|
\\
&&{}+\sup_{u\in\mathbb{R}}\bigl|\gamma_k(u)\bigr|= O_P(1)
\sqrt{k} n^{-1/4+\kappa
/2}+o(1)k^{-(m-2)/2-1/12}\\
&=&o_P(1),
\end{eqnarray*}
where the facts that $\sup_{u\in\mathbb{R}}\|Z_k (u)\|\le C\sqrt{k}$
and $\sup_{u\in\mathbb{R}}|\gamma_k(u)|\le\break  C k^{-(m-2)/2-1/12}$ with
some constant $C>0$ are given in Lemma C.1 in the supplementary
material of the paper.

For the normality, in view of the consistency of $\widehat{\sigma}_e$
and $\widehat{L}_{n1}(1,0)$, we show the result with the replacement of
$\sigma_e$ and $L_1(1,0)$. Meanwhile, in order to correspond to the
plug-in of $\widehat{\theta}_n$, denote by $\widehat{Z}$ the matrix $Z$
defined in Section~\ref{sec2} with replacement of $\theta_0$ by $\widehat{\theta
}_n$. Noting that $\widetilde{c}= \widetilde{c} (\theta_0)=(Z^\top
Z)^{-1}Z^\top Y$ and $Y=Zc+\gamma+e$ given in Section~\ref{sec2.1},
\begin{eqnarray*}
\widehat{c}&=&\widetilde{c}(\widehat{\theta}_n)=\bigl(
\widehat{Z}^\top \widehat{Z}\bigr)^{-1}\widehat{Z}^\top(Zc+
\gamma+e)
\\
&=& c+\bigl(\widehat{Z}^\top\widehat{Z}\bigr)^{-1}
\widehat{Z}^\top(\gamma +e)+\bigl(\widehat{Z}^\top\widehat{Z}
\bigr)^{-1}\widehat{Z}^\top(Z-\widehat{Z})c.
\end{eqnarray*}
It follows from Lemma~\ref{lemma2} that
\begin{eqnarray*}
&&\widehat{g}_n(u)-g(u)\\
&&\qquad=Z_k^\top(u)
\widehat{c}-g(u)=Z_k^\top (u) (\widehat{c}-c) -
\gamma_k(u)
\\
&&\qquad=Z_k^\top(u) \bigl(\widehat{Z}^\top\widehat{Z}
\bigr)^{-1}\widehat{Z}^\top (\gamma+e)+Z_k^\top(u)
\bigl(\widehat{Z}^\top\widehat{Z}\bigr)^{-1}\widehat
{Z}^\top(Z-\widehat{Z})c -\gamma_k(u)
\\
&&\qquad=\frac{1}{\sqrt{n}}L^{-1}_1(1,0)Z_k^\top(u)
\widehat{Z}^\top e\bigl(1+o_P(1)\bigr)+\frac{1}{\sqrt{n}}
L^{-1}_1(1,0) Z_k^\top(u)\widehat
{Z}^\top\gamma
\\
&&\qquad\quad{}+\frac{1}{\sqrt{n}} L^{-1}_1(1,0) Z_k^\top(u)
\widehat{Z}^\top (Z-\widehat{Z})c -\gamma_k(u)
\\
&&\qquad=\frac{1}{\sqrt{n}L_1(1,0)}Z_k^\top(u)Z^\top e+
\frac{1}{\sqrt
{n}L_1(1,0)}Z_k^\top(u) (\widehat{Z}
-Z)^\top e
\\
&&\qquad\quad{}+\frac{1}{\sqrt{n}L_1(1,0)} Z_k^\top(u)\widehat{Z}^\top
\gamma\\
&&\qquad\quad{} +\frac{1}{\sqrt{n}L_1(1,0)} Z_k^\top(u)\widehat{Z}^\top(Z-
\widehat {Z})c -\gamma_k(u).
\end{eqnarray*}

To fulfill the normality, we need to show
\begin{eqnarray*}
&&(1)\quad \sigma_e^{-1} L_1(1,0)^{-1/2}
\frac{1}{\sqrt[4]{n}}\bigl\|Z_k(u)\bigr\| ^{-1} Z_k^\top(u)Z^\top
e\to_D N(0,1),
\\
&&(2)\quad\frac{1}{\sqrt[4]{n}\|Z_k(u)\|}Z_k^\top(u)\widehat{Z}^\top
\gamma=o_P(1),\\
&& (3)\quad \frac{1}{\sqrt[4]{n}\|Z_k(u)\|}Z_k^\top
(u)\widehat{Z}^\top(Z-\widehat{Z})c=o_P(1),
\\
&&(4)\quad \sqrt[4]{n}\bigl\|Z_k(u)\bigr\|^{-1}\gamma_k(u)=o(1),\\
&&
(5)\quad\frac
{1}{\sqrt[4]{n}}\bigl\|Z_k(u)\bigr\|^{-1}Z_k^\top(u)
(\widehat{Z} -Z)^\top e=o_P(1).
\end{eqnarray*}

To begin with (1), observe that
\begin{eqnarray*}
&&n^{-1/4}\sigma_e^{-1}L^{-1/2}_1(1,0)
\bigl\|Z_k(u)\bigr\|^{-1}Z_k^\top(u)Z^\top
e
\\
&&\qquad= n^{-1/4}\sigma_e^{-1}L^{-1/2}_1(1,0)
\bigl\|Z_k(u)\bigr\|^{-1} \sum_{t=1}^n
Z_k^\top(u) Z_k\bigl(\theta_0^\top
x_t\bigr)e_t,
\end{eqnarray*}
which is a martingale array in view of Assumption \ref{assa}. We shall use
Corollary~3.1 of \citet{peterhall1980} to show the normality of (1).

The conditional variance process is, by Lemma~\ref{lemma2},
\begin{eqnarray*}
&&\frac{1}{\sqrt{n}}\sigma_e^{-2}L^{-1}_1(1,0)
\bigl\|Z_k(u)\bigr\|^{-2} \sum_{t=1}^n
\bigl(Z_k^\top(u) Z_k\bigl(
\theta_0^\top x_t\bigr)\bigr)^2 E
\bigl(e_t^2|\mathcal {F}_{n,t-1}\bigr)
\\
&&\qquad=\frac{1}{\sqrt{n}}L^{-1}_1(1,0)\bigl\|Z_k(u)
\bigr\|^{-2} \sum_{t=1}^n
\bigl(Z_k^\top (u) Z_k( x_{1t})
\bigr)^2
\\
&&\qquad=\frac{1}{\sqrt{n}}L^{-1}_1(1,0)\bigl\|Z_k(u)
\bigr\|^{-2} Z_k^\top(u) \Biggl(\sum
_{t=1}^n Z_k( x_{1t})^\top
Z_k(x_{1t}) \Biggr)Z_k (u)
\\
&&\qquad=L^{-1}_1(1,0)\bigl\|Z_k(u)\bigr\|^{-2}
Z_k^\top(u) \biggl(\frac{1}{\sqrt
{n}}Z^\top Z
\biggr) Z_k(u)
\\
&&\qquad=\bigl\|Z_k(u)\bigr\|^{-2} Z_k^\top(u)
Z_k(u) \bigl(1+o_P(1)\bigr)=1+o_P(1).
\end{eqnarray*}
Moreover, to make the conditional Lindeberg's condition fulfilled it
suffices to show
\begin{eqnarray*}
&&\bigl\|Z_k(u)\bigr\|^{-4}\frac{1}{n}\sum
_{t=1}^n E\bigl[\bigl(Z_k^\top(u)Z_k(x_{1t})
e_t\bigr)^4| \mathcal{F}_{n,t-1}\bigr]
\\
&&\qquad\le C\bigl\|Z_k(u)\bigr\|^{-4}\frac{1}{n}\sum
_{t=1}^n \bigl\|Z_k(u)\bigr\|^4
\bigl\|Z_k(x_{1t})\bigr\| ^4=C \frac{1}{n} \sum
_{t=1}^n \bigl\|Z_k(x_{1t})
\bigr\|^4=o_P(1)
\end{eqnarray*}
by a routine calculation using the density of $t^{-1/2}x_{1t}$ in Lemma~\ref{lemma5a}. This finishes the normality for (1). For the sake of
brevity, the proof for (2)--(5) is relegated to the supplementary
material. The outline then is completed.
\end{pf*}

\begin{pf*}{Outline of the Proof of Theorem~\ref{th4}}
Denote for any $\vartheta=(\beta,\theta)$,
\begin{eqnarray*}
\mathfrak{S}_n(\vartheta)&=& %
\pmatrix{\mathfrak{S}_{n,1}(
\vartheta)\vspace*{2pt}
\cr
\mathfrak {S}_{n,2}(\vartheta) } %
=
\pmatrix{\displaystyle\frac{\partial L_n ( \vartheta)}{\partial\beta}
\vspace*{2pt}\cr
\displaystyle\frac
{\partial L_n(\vartheta)}{\partial\theta} } %
,
\\
\mathfrak{J}_n(\vartheta)&=& %
\pmatrix{\mathfrak{J}_{n,11}
(\vartheta) & \mathfrak {J}_{n,12}(\vartheta)\vspace*{2pt}\cr
\mathfrak{J}_{n,21}(\vartheta) & \mathfrak {J}_{n,22}(\vartheta)
} %
= %
\pmatrix{\displaystyle\frac{\partial^2 L_n(\vartheta)}{\partial\beta\,
\partial\beta^\top}&\displaystyle\frac{\partial^2L_n(\vartheta)}{\partial\beta\,
\partial\theta^\top}
\vspace*{2pt}
\cr
\displaystyle\frac{\partial^2 L_n(\vartheta) }{\partial
\theta\,\partial\beta^\top}&\displaystyle\frac{\partial^2 L_n(\vartheta)}{\partial
\theta\,\partial\theta^\top} } %
.
\end{eqnarray*}
Also, for any $\mu=(\lambda,\alpha)$, $\mathfrak{S}_n(\mu)$ and
$\mathfrak{J}_n (\mu)$ are defined similarly but with the parameters rotated.

Denote $\widetilde{D}_n=\mathrm{diag}(nI_d, D_n)$. Thus, \eqref{c2b} may
be equivalently written as
%
%eA.2 #&#
\begin{equation}
\label{th4b1} \widetilde{D}_n^{-1}\mathfrak{S}_n(
\mu_0)+\widetilde{D}_n^{-1}\mathfrak
{J}_n(\mu_n) \widetilde{D}_n^{-1}
\widetilde{D}_n(\widehat{\mu}_n -\mu_0)=0.
\end{equation}

It follows from \eqref{th4b1} that
%
%eA.3 #&#
%eA.4 #&#
\begin{eqnarray}\label{th4c}
&&n^{-1}\mathfrak{S}_{n,1}(\mu_0)+n^{-2}
\mathfrak{J}_{n,11}(\mu _n)n(\widehat{\lambda}_n
-\lambda_0)
\nonumber
\\[-8pt]
\\[-8pt]
\nonumber
&&\qquad{} + n^{-1}\mathfrak{J}_{n,12}(\mu_n)D_n^{-1}D_n
(\widehat{\alpha }_n -\alpha_0)=0,
\\
\label{th4d}
&&D_n^{-1}\mathfrak{S}_{n,2}(
\mu_0)+D_n^{-1}\mathfrak{J}_{n,21}(
\mu_n)n^{-1} n(\widehat{\lambda}_n -
\lambda_0)
\nonumber
\\[-8pt]
\\[-8pt]
\nonumber
&&\qquad +D_n^{-1}\mathfrak{J}_{n,22}(
\mu_n)D_n^{-1}D_n( \widehat{\alpha
}_n -\alpha_0)=0.
\end{eqnarray}

The results of \eqref{th4a} and \eqref{th4b} will be derived from \eqref
{th4c} and \eqref{th4d}, respectively. These are shown in the following
two steps.

\textit{Step} 1: We first prove \eqref{th4b} from \eqref{th4d}. First of
all, note that $\mathfrak{S}_{n,2}(\mu_0)$ and $\mathfrak{J}_{n,22}(\mu
_0)$ are exactly the $S_n(\alpha_0)$ and $J_n(\alpha_0)$ in Theorem~\ref
{th31}, respectively, since $y_t-\beta_0^\top x_t$ in model \eqref{a1}
plays the same role as $y_t$ in model \eqref{a2}. Therefore,
%
%eA.5 #&#
\begin{equation}
\label{th4d1} D_n^{-1}\mathfrak{S}_{n,2}(
\mu_0)\to_D R^{1/2}W(1) \quad\mbox {and}\quad
D_n^{-1}\mathfrak{J}_{n,22}(\mu_0)D_n^{-1}
\to_P R,
\end{equation}
where $R$ and $W$ are defined in Theorem~\ref{th31}.

To prove \eqref{th4b}, it therefore suffices to show that
%
%eA.6 #&#
\begin{equation}
\label{th4e} D_n( \widehat{\alpha}_n -
\alpha_0)=\bigl[D_n^{-1}\mathfrak{J}_{n,22}(
\mu_0) D_n^{-1}\bigr]^{-1}D_n^{-1}
\mathfrak{S}_{n,2}(\mu_0)+o_P(1).
\end{equation}

Once again, Theorem~10.1 of \citet{wooldridge1994} is used to prove
\eqref{th4e}. It is shown in detail in the supplemental material that
$n( \widehat{\lambda}_n -\lambda_0)= Q^\top n(\widehat{\beta}_n -\beta
_0 )=O_P(1)$ and $D_n^{-1} \mathfrak{J}_{n,21}(\mu_0) n^{-1}=o_P(1)$.
Define for some $\delta>0$, $\widetilde{C}_n= n^{-\delta} \widetilde
{D}_n=\operatorname{diag}( n^{1-\delta}I_d, C_n)$ such that $\widetilde{C}_n
\widetilde{D}_n^{-1}\to0$ as $n\to\infty$, where $C_n= n^{-\delta}
D_n$ used in the proof of Theorem~\ref{th31}. It follows from \eqref
{th4d} that
\begin{eqnarray*}
0&=&D_n^{-1}\mathfrak{S}_{n,2}(
\mu_0)+n^{-2\delta}C_n^{-1}\bigl[\mathfrak
{J}_{n,21}(\mu_n)- \mathfrak{J}_{n,21}(
\mu_0)\bigr]n^{-1+\delta}n(\widehat {\lambda}_n -
\lambda_0)+o_P(1)
\\
&&{}+D_n^{-1}\mathfrak{J}_{n,22}(\mu
_0)D_n^{-1}D_n( \widehat{
\alpha}_n -\alpha_0)
\\
&&{}+n^{-2\delta}C_n^{-1}\bigl[\mathfrak{J}_{n,22}(
\mu_n)-\mathfrak {J}_{n,22}(\mu_0)
\bigr]C_n^{-1} D_n( \widehat{
\alpha}_n -\alpha_0).
\end{eqnarray*}

The requirements (i)--(ii) in Theorem~10.1 of \citet{wooldridge1994} are
trivially fulfilled and the requirement (iii) will be satisfied if we
can show
%
%eA.7 #&#
%eA.8 #&#
\begin{eqnarray}
\sup_{\{\mu: \|\widetilde{C}_n(\mu-\mu_0)\|<1\}} \bigl\|n^{-1+\delta} C_n^{-1}
\bigl[\mathfrak{J}_{n,21}(\mu) -\mathfrak{J}_{n,21} (
\mu_0)\bigr]\bigr\| &=&o_P(1),\label{th4f}
\\
\sup_{\{\mu: \|\widetilde{C}_n(\mu-\mu_0)\|<1\}} \bigl\|C_n^{-1} \bigl[\mathfrak
{J}_{n,22} (\mu)-\mathfrak{J}_{n,22}(\mu_0)\bigr]
C_n^{-1}\bigr\| &=&o_P(1).\label{th4g}
\end{eqnarray}
With the choice of $\delta\in(0,1/24)$, both \eqref{th4f} and \eqref
{th4g} are proved in the supplemental material of the paper, and hence
the requirement \eqref{th4e} is verified if we choose $\delta\in(0, 1/24)$.

Furthermore, \eqref{th4d1} shows that condition (iv) in Wooldridge's
theorem holds. Thus, the limit distribution \eqref{th4b} now follows directly.

\textit{Step} 2: We now turn to prove \eqref{th4a} from \eqref{th4c}.
Since $\mathfrak{J}_{n,12} (\mu_n) =\mathfrak{J}_{n,21}(\mu_n)^\top$,
$D_n( \widehat{\alpha}_n -\alpha_0)=O_P(1)$ by Step 1,
$D_n^{-1}\mathfrak{J}_{n,21}(\mu_0)n^{-1}=o_P(1)$ (as shown in the
supplementary material), and $\mathfrak{J}_{n,11}(\mu_n)$ is
independent of $\mu_n$, we have
\[
n(\widehat{\lambda}_n -\lambda_0) =\bigl(n^{-2}
\mathfrak{J}_{n,11}(\mu _n)\bigr)^{-1}
n^{-1} \mathfrak{S}_{n,1} (\mu_0)+o_P(1)
\]
by \eqref{th4f}. Note that
\begin{eqnarray*}
\frac{1}{n}\mathfrak{S}_{n,1}(\mu_0)&=&
\frac{1}{n}Q^\top\frac{\partial
L_n(\vartheta_0)}{\partial\beta}=\frac{1}{n}Q^\top
\sum_{t=1}^n\bigl(y_t-
\beta_0^\top x_t-\widehat{g}_n
\bigl(\theta_0^\top x_t\bigr)
\bigr)x_t
\\
&=&\frac{1}{n} \sum_{t=1}^ne_tQ^\top
x_t+o_P(1) \to_D Q^\top\int
_0^1V(r)\,dU(r),
\end{eqnarray*}
as shown in Appendix D of the supplementary material when $n\to\infty
$. Note also that
\[
\frac{1}{n^2} \mathfrak{J}_{n,11}(\mu_n)=
\frac{1}{n^2}Q^\top\sum_{t=1}^n
x_t x_t^\top Q\to Q^\top\int
V(r)V(r)^\top \,dr Q
\]
almost surely using Theorem~3.1 of \citet{phillips2001}, from which
\eqref{th4a} follows. The outline of the proof is completed.
\end{pf*}

\begin{pf*}{Proof of Theorem~\ref{th5}}
The result of \eqref{th5a} follows directly from \eqref{th4a}. In view
of the proof of Theorem~\ref{th32} as well as \eqref{th4b}, equation
\eqref{th5b} holds.
\end{pf*}

\begin{pf*}{Proof of Theorem~\ref{th6}}
In view of \eqref{th4b} and the proofs of Theorems \ref{th32}--\ref
{th33}, it holds.
\end{pf*}
\end{appendix}

%sA #&#
\section*{Acknowledgements}
The authors are grateful to the Editor, Professor Runze Li, an
Associate Editor and three referees for their valuable and constructive
comments and suggestions that substantially improve earlier versions of
the paper. Thanks also go to many seminar and conference participants
for their useful comments when this work was presented. The authors
acknowledge the financial support from the Norwegian Research Council.

\begin{supplement}[id=suppA]
\stitle{Additional technical details}
\slink[doi]{10.1214/15-AOS1372SUPP} %[doi,text={...}] - jei reikia suskaldyti doi
\sdatatype{.pdf}
\sfilename{aos1372\_supp.pdf}
\sdescription{The proofs and technical details
that are omitted in the paper are provided in the supplement that
accompanies the paper.}
\end{supplement}

%\bibliography{123}
%
% imsref loaded by akundreckaite, 2015-10-06 14:37:19

%\begin{appendix}
%\section{}
%\end{appendix}

%\section*{Acknowledgments}

%\begin{thebibliography}{99}
%\bibitem[\protect\citeauthoryear{}{}]{r1}
%\bibitem{r1}
%\end{thebibliography}

\printaddresses
\end{document}